\documentclass[full]{siamart171218}

\pdfoutput=1

\usepackage[section]{placeins} 

\usepackage{algorithmic}
\usepackage[caption=false]{subfig}
\ifpdf
  \DeclareGraphicsExtensions{.eps,.pdf,.png,.jpg}
\else
  \DeclareGraphicsExtensions{.eps}
\fi

\usepackage{amssymb}
\usepackage{amsmath}
\usepackage{amsfonts}
\usepackage{upgreek}
\usepackage{siunitx}
\usepackage{listings}
\usepackage{multirow}
\usepackage{multicol}
\usepackage{bigdelim}
\usepackage{color}
\usepackage{stmaryrd}
\usepackage{soul}
\usepackage[symbol]{footmisc}
\usepackage{booktabs}
\usepackage[section]{placeins} 
\usepackage[titletoc]{appendix}
\usepackage{calc}
\usepackage[capitalize]{cleveref}

\usepackage{graphicx}
\usepackage{graphics}
\usepackage{float}
\usepackage{wrapfig}
\usepackage[percent]{overpic}
\usepackage{varwidth}
\usepackage{epstopdf}
\usepackage{adjustbox}
\usepackage{tikz}
\usepackage{pgfplots}
\usetikzlibrary{arrows,matrix,positioning,fit}
\usetikzlibrary{shapes,positioning}
\usetikzlibrary{backgrounds}
\pgfplotsset{compat=1.14}
\usepgfplotslibrary{colorbrewer}
\usepgfplotslibrary{patchplots}
\usepgfplotslibrary[colorbrewer]
\usetikzlibrary{pgfplots.colorbrewer}
\usetikzlibrary[pgfplots.colorbrewer]
\usepgfplotslibrary{fillbetween}
\usepgfplotslibrary{units}
\usetikzlibrary{spy}
\usepackage{pgfplotstable}
\graphicspath{ {Figs/} }
\usetikzlibrary{external}

\newcommand{\etal}{et~al.}

\newcommand{\hfd}{\hat{f}^{\delta}}

\newcommand{\hfI}{\hat{f}^{\delta I}}
\newcommand{\hfC}{\hat{f}^{\delta C}}

\newcommand{\hud}{\hat{u}^{\delta}}

\newcommand\px[2]{\frac{\partial #1}{\partial {#2}}}
\newcommand\pxi[3]{\frac{\partial^{#1}#2}{\partial {#3}^{#1}}}
\newcommand\pxvar[2]{\partial_{#2} #1}

\newlength\myheight
\newlength\mydepth
\settototalheight\myheight{Xygp}
\newsiamremark{remark}{Remark}
\newsiamremark{hypothesis}{Hypothesis}
\crefname{hypothesis}{Hypothesis}{Hypotheses}
\newsiamthm{claim}{Claim}

\headers{\emph{p}-multigrid Fourier Analysis}{W. Trojak and F.D. Witherden}

\title{On Fourier analysis of polynomial multigrid for arbitrary multi-stage cycles}

\author{W. Trojak\thanks{Department of Ocean Engineering, Texas A\&M University, College Station, TX 77843 (\email{wt247@tamu.edu}, \email{fdw@tamu.edu}).}
\and F.D. Witherden\footnotemark[2]}

\begin{document}

\maketitle

\begin{abstract}
    The Fourier analysis of the \emph{p}-multigrid acceleration technique is considered for a dual-time scheme applied to the advection-diffusion equation with various cycle configurations. It is found that improved convergence can be achieved through \emph{V}-cycle asymmetry where additional prolongation smoothing is applied. Experiments conducted on the artificial compressibility formulation of the Navier--Stokes equations found that these analytic findings could be observed numerically in the pressure residual, whereas velocity terms---which are more hyperbolic in character---benefited primarily from increased pseudo-time steps.
\end{abstract}

\begin{keywords}
High-order, Flux reconstruction, Multigrid, Dual-time, Fourier analysis 
\end{keywords}

\begin{AMS}
65M60, 65T99, 65M55, 76D99
\end{AMS}



\section{Introduction}\label{sec:intro}
The artificial compressibility method~(ACM)~\cite{Chorin1967} is a means of solving the incompressible Navier--Stokes equations in a manner that is compatible with compressible solvers. The most widely applied method for the incompressible Navier--Stokes equations is the pressure correction method where pressure corrections from a Poisson equation are propagated into the weakly coupled velocity field. This method has the disadvantage of indirect communication which can reduce parallel efficiency. ACM instead couples pressure to the continuity equation and consequently has seen increasing popularity for computational fluid dynamics; however, for each time step it does require that the artificial pressure waves are allowed to propagate in pseudo-time such that the converged, incompressible solution is reached. A technique commonly used to achieve this converged state is dual-time stepping~\cite{Chorin1967,Peyret1976}, and due to the requirement to converge the system for each time step, it follows that an implicit temporal integration scheme is applied. Other approaches have been explored, such as solving the linearised pseudo-time system with GMRES~\cite{Rogers1995}. However, this method requires preconditioning and has parallelisation issues common with these implicit methods.

Relative to pseudo-time, the system is driven to a steady state, and hence many convergence acceleration techniques are applicable. Several approaches have been developed, notably simple spatially-varying time steps, alternating direction implicit schemes~\cite{Peaceman1955}, implicit-explicit hybrid schemes~\cite{Hsu2002}, and the use of complex relaxation schemes such as LU-SSOR~\cite{Yoon1987}. The technique which is the concern of the present work is the multigrid method~\cite{Arnone1993} which is particularly effective for elliptic problems and hence may be well suited to accelerating ACM due to the nature of the artificial pressure waves.

Important to the application of multigrid acceleration is which spatial scheme is employed. We are interested in the use of spectral element discretisations and, in particular, the flux reconstruction method~(FR)~\cite{Huynh2007} which can be understood as a generalisation of the nodal discontinuous Galerkin approach~\cite{Hesthaven2008}. This method is of interest due to its high-order and globally unstructured nature combined with locally structured compute that lends itself to  modern computer architectures~\cite{Witherden2014}. High-order methods are particularly beneficial in the context of ACM due to the lack of solution discontinuities, hence making these techniques highly efficient in the approximation in spatial derivatives. 

The application of multigrid methods---such as geometric multigrid---is complicated by the unstructured formulation of FR. However, the high spatial order lends itself readily to \emph{p}-multigrid acceleration methods where for the same element coarser levels are introduced via restricting the solution to lower polynomial orders. There is a rich body of literature considering the Fourier analysis of geometric multigrid methods, with analysis advancing to more general deep cycles such as the work of Wienands~\etal~\cite{Wienands2001}, where it was theoretically shown that contraction factors could deteriorate for schemes with more stages due to aliasing on the coarsest levels. We wish to develop a theoretical framework to explore the effect of cycle design on acceleration of \emph{p}-multigrid methods.
\section{The FR Approach}\label{sec:fr}
    The analysis of the methods to be presented will at times require the explicit coupling of temporal integration methods to a spatial scheme to produce the eigenvalues of the system. The spatial scheme used is the  FR~\cite{Huynh2007,Vincent2010} method which lies within the set of discontinuous spectral element methods. For the purpose of this analysis, the FR method is used for approximating the first derivative of a function, with second derivatives handled through the introduction of auxiliary variables. Let us set the function $f$ such that $f(u): \mathbb{R} \mapsto \mathbb{R}$, and the domain of the spatial variable $x\in\Omega$. The spatial domain is subdivided into sub-domains $\Omega_i$, such that $\bigcup^N \Omega_i = \Omega$ and $\Omega_i \cap \Omega_j = \emptyset$ if $i \neq j$. In one dimension, we define a reference element and variable, $\xi \in \hat{\Omega} = [-1,1]$, for which we introduce the Jacobian $J_j: \Omega_j \mapsto \hat{\Omega}$.

    If we have a solution $u(x)$ and function $f(u)$, FR forms a degree $p$ polynomial approximation of $f$ in $\Omega_i$ transformed to the reference space via the values at a set of $p+1$ nodal points $\xi_{j}$. We denote the discontinuous approximation as $\hfd_i$
    \[
        \hfd_i = \sum^p_{j=0} f\big(\hud_i(\xi_j)\big)l_j(\xi) \quad \mathrm{where} \quad l_j = \prod^{p}_{\substack{k=0 \\k\neq j}}\frac{\xi-\xi_k}{\xi_j-\xi_k},
    \]
    which is similarly defined for $\hud_i$. The FR methodology is then concerned with updating this polynomial such that the approximation is  $C^0$ continuous between elements. This is achieved via
    \[
        \hfC = \hfd + (\hfI_{i,L} - \hfd_{i,L})h_L + (\hfI_{i,R} - \hfd_{i,R})h_R,
    \]
    where $\hfd_{i,L} = \hfd_i(-1)$ is the interpolated value of $\hfd_i$ at the left interface, and $\hfI_{i,L}$ is the common interface function value at the left interface. Similar definition follow for the right interface. For hyperbolic problems, the interface flux may be found by using information from the adjacent cell to pose a Riemann problem. There are many appropriate methods for the approximation or solution of these problems~\cite{Toro2009}, and it has also been demonstrated~\cite{Jameson2011} that the \emph{E-flux} condition is important in the proof of stability.  The functions $h_L$ and $h_R$ are correction functions with the boundary conditions $h_L(-1) = h_R(1) = 1$ and $h_L(1) = h_R(-1) = 0$, and if they are set to left and right Radau polynomials then a nodal DG scheme is recovered~\cite{Huynh2007}. With this, the spatial derivative can straightforwardly be obtained, and if $h_L \in \mathbb{P}^{p+1}$, then it is possible for $\partial\hfC_i/\partial\xi \in \mathbb{P}^p$. 
    
\section{Pseudo-Time Stepping}\label{sec:pseudo}
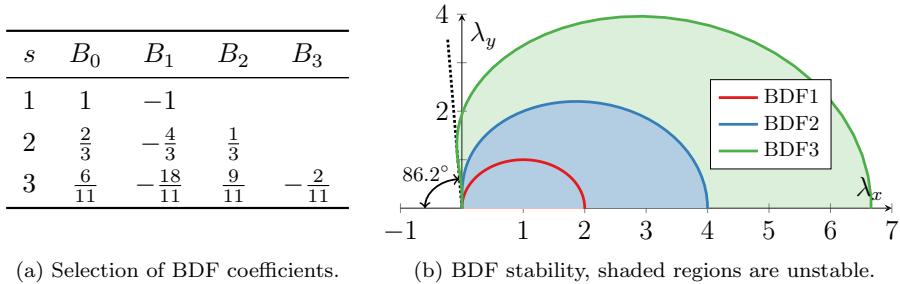
\begin{figure}[tbhp]
        \centering
        \subfloat[Selection of BDF coefficients.]{\label{tab:bdf}
            \adjustbox{width=0.35\linewidth,valign=b}{\begin{tabular}{c c c c c}
		        \toprule
		        $s$ & $B_0$ & $B_1$ & $B_2$ & $B_3$ \\ \midrule
		        1 & $1$ & $-1$ &  &  \\[0.75ex]
		        2 & $\frac{2}{3}$ & $-\frac{4}{3}$ & $\frac{1}{3}$ &  \\[0.75ex]
		        3 & $\frac{6}{11}$ & $-\frac{18}{11}$ & $\frac{9}{11}$ & $-\frac{2}{11}$ \\
 		        \bottomrule\vspace{0.09cm}
            \end{tabular}}
        }
        ~
        \subfloat[BDF stability, shaded regions are unstable.]{\label{fig:bdf_base_stab}\adjustbox{width=0.55\linewidth,valign=b}{\begin{tikzpicture}
	\begin{axis}[name=plot1,xlabel={$\lambda_x$},ylabel={$\lambda_y$},
	    axis line style={latex-latex},
        axis y line=middle,
        axis x line=middle,
        width=0.7\textwidth,
        height=0.35\textwidth,
        xmin=-1,
	    xmax=7,
        ymin=0,
    	ymax=4,
    	xtick={-1,0,1,2,3,4,5,6,7},
    	ytick={0,1,2,3,4},
    	yticklabels={ , ,2, ,4},
    	legend style={at={(0.63,0.2)},anchor=south west,font=\small},
    	legend cell align={left},
    	style={font=\large}
    	]
    	\path[name path=axis] (axis cs:-1,0) -- (axis cs:7,0);
    	
        \addplot[color={black}, style={very thick,densely dotted}, forget plot] coordinates{(0,0) (3.5*tan(86.2-90),3.5)};
        \node[style={font=\small}] at (axis cs: -0.58,0.72) {$86.2^\circ$};
        \draw[thick, <->] (axis cs:-0.6,0) arc(180:93.8:0.6);
    	
        \addplot[name path=BDF1,color={Set1-A}, style={very thick},line legend]table[x=x,y=y,col sep=comma,unbounded coords=jump]{./Figs/data/BDF1_region.csv};
        
        \addplot[name path=BDF2,color={Set1-B}, style={very thick},line legend]table[x=x,y=y,col sep=comma,unbounded coords=jump]{./Figs/data/BDF2_region.csv};
        
        \addplot[name path=BDF3,color={Set1-C}, style={very thick},line legend]table[x=x,y=y,col sep=comma,unbounded coords=jump]{./Figs/data/BDF3_region.csv};
        
        \addplot[very thick,color={Set1-C},fill={Set1-C},fill opacity=0.2,forget plot] fill between[of=BDF3 and BDF2,soft clip={domain=-1:7}];
        
 		\addplot[very thick,color={Set1-B},fill={Pastel1-B},fill opacity=1,forget plot] fill between[of=BDF2 and BDF1,soft clip={domain=-1:7}];
		
 		\addplot[very thick,color={Set1-A},fill={Pastel1-A},fill opacity=1,forget plot] fill between[of=BDF1 and axis,soft clip={domain=-1:7}];
		
        \legend{BDF1,BDF2,BDF3}
    \end{axis}    		
\end{tikzpicture}}}
        \caption{\label{fig:bdf}BDF Schemes.}
    \end{figure}

    To introduce the dual-time method, consider the ordinary differential equation~(ODE)
    \begin{equation}\label{eq:linear_ode}
        \px{u}{t} -\lambda u = 0 \quad \mathrm{for} \quad (x,t)\in \Omega\times\mathbb{R}_+,
    \end{equation}
    which may be modified to incorporate pseudo-time terms as
	\begin{equation}\label{eq:linear_pseduo}
		\px{u}{\tau} + \px{u}{t} - \lambda u = 0 \quad \mathrm{for} \quad (x,t,\tau)\in \Omega\times\mathbb{R}_+^2,
	\end{equation}
	such that when a steady state in pseudo-time is reached, then a solution to \cref{eq:linear_ode} is reached. To simplify later analysis, we will restrict the spatial domain to be periodic, thus restricting the equation to an initial value problem. To solve this system, we will employ explicit Runge--Kutta (ERK) integration in pseudo-time.  Such schemes may be defined through a Butcher tableau~\cite{Butcher1964} as
	\begin{equation}\label{eq:butcher}
	    \arraycolsep=2.5pt\def\arraystretch{1.2}
		\begin{array}{c|c}
			\mathbf{c} & \mathbf{A} \\ \hline
			& {\mathbf{b}^T}
		\end{array}
	\end{equation}
    For ERK schemes, the coefficient matrix $\mathbf{A}$ is strictly lower triangular. The ERK scheme applied to integration of the ODE in \cref{eq:linear_ode} can be written as
    \begin{equation}
            u_{n+1} = u_{n} + \sum^r_{i=1}\Delta t q_i, \quad \mathrm{with} \quad q_i = \lambda\bigg(u_{n} + \Delta t\sum^{i-1}_{j=1}a_{ij}q_j\bigg),
    \end{equation}
    where $\Delta t$ is the time step size. For the system presented in \cref{eq:linear_pseduo}, this ERK scheme will be used for the pseudo-time integration, whereas physical time stepping will be performed with the implicit backward-difference formulae (BDF). The general form for a degree $s$ BDF scheme can be expressed as
    \begin{equation}
        u_{n+1} = -\sum^{s-1}_{i=0}B_{i+1}u_{n-i} + \Delta tB_0\lambda u_{n+1}.
    \end{equation}
    Example coefficients and stability regions for several BDF schemes are available in \cref{fig:bdf}. 
    
    The implicit and explicit integrators for physical- and pseudo-time can now be combined to calculate the solution advanced by the pseudo step, $\Delta\tau$, thus giving the following system of equations
	\begin{subequations}
		\begin{align}
			u_{n+1, m+1} =&\: u_{n+1, m} + \sum^{r}_{i=1}\frac{\Delta\tau b_i}{\alpha_{PI}}q_i, \\
			q_i =&\: \lambda \Big(u_{n+1,m} + \Delta\tau\sum^{i-1}_{j=1}a_{i,j}q_j\Big) - \frac{1}{\Delta tB_0}\Big(u_{n+1,m} + \sum^{s-1}_{l=0}B_{l+1}u_{n-l}\Big).
		\end{align}
    \end{subequations}
    From the use of explicit pseudo-time stepping, it logically follows that we assume $\Delta\tau \ll \Delta t$, and hence the term $\alpha_{PI}=1+b_i\Delta\tau/B_0\Delta t\rightarrow1$ can be neglected. We now wish manipulate this into a matrix form to facilitate our later work; applying the terms of \cref{eq:butcher}, the following is obtained
	\begin{equation*}
		\mathbf{q} = \lambda u_{n+1,m}\mathbf{e} + \lambda\Delta\tau\mathbf{A}\mathbf{q} - \frac{1}{\Delta tB_0}\Big(u_{n+1,m} + \sum^{s-1}_{l=0}B_{l+1}u_{n-l}\Big)\mathbf{e},
	\end{equation*}
    where $\mathbf{q} = [q_1,\dots,q_r]^T$, and $\mathbf{e} = [1,\dots,1]^T$. This, in turn, implies
	\begin{subequations}
		\begin{align*}
			u_{n+1,m+1} =&\: u_{n+1,m} + \Delta\tau\mathbf{b}^T\mathbf{q}, \\
			\mathbf{q} =&\: (\mathbf{I} - \lambda\Delta\tau\mathbf{A})^{-1}\bigg[\lambda u_{n+1,m} - \frac{1}{\Delta tB_0}\Big(u_{n+1,m} + \sum^{s-1}_{l=0}B_{l+1}u_{n-l}\Big)\bigg]\mathbf{e}.
		\end{align*}
	\end{subequations}
    To obtain the system amplification factor, we factorise $\mathbf{q}$ in terms of $u_{n+1,m}$ by initially separating the pseudo-time amplification and source terms as
	\begin{equation*}
		\mathbf{q} = (\mathbf{I} - \lambda\Delta\tau\mathbf{A})^{-1}\mathbf{e}\bigg(\lambda - \frac{1}{\Delta tB_0}\bigg)u_{n+1,m}  -\frac{(\mathbf{I} - \lambda\Delta\tau\mathbf{A})^{-1}\mathbf{e}}{\Delta tB_0}\sum^{s-1}_{l=0}B_{l+1}u_{n-l}.
	\end{equation*}
	Therefore,
	\begin{multline}\label{eq:ps1}
		u_{n+1,m+1} = \underbrace{\Bigg[1 + \bigg(\lambda\Delta\tau - \frac{\Delta\tau}{\Delta tB_0}\bigg)\mathbf{b}^T(\mathbf{I} - \lambda\Delta\tau\mathbf{A})^{-1}\mathbf{e}\Bigg]}_Pu_{n+1,m}, \\
		- \underbrace{\frac{\Delta\tau\mathbf{b}^T(\mathbf{I} - \lambda\Delta\tau\mathbf{A})^{-1}\mathbf{e}}{\Delta tB_0}}_C\sum^{s-1}_{l=0}B_{l+1}u_{n-l},
	\end{multline}	
	and it can be seen that this is purely a function of the ERK and BDF schemes, together with the factors $\lambda\Delta\tau$ and $\Delta\tau/\Delta t$.
	
	To demonstrate the effect of the coupled system, we present the stability regions of \cref{eq:ps1} as the pseudo step number, $m$, is varied. This was calculated using the amplification factor, defined as
	\begin{subequations}
	    \begin{align}
	        \frac{u_{n+1,M}}{u_{n}} &= P^M - \bigg[\sum^{M-1}_{j=0}P^j\bigg]C\sum^{s-1}_{l=0}B_{l+1}\exp{(\lambda l\Delta t)},\\
	        &= P^M - \frac{1-P^M}{1-P}C\sum^{s-1}_{l=0}B_{l+1}\exp{(\lambda l\Delta t)}.
	    \end{align}
	\end{subequations}
    In the second step we have assumed $|P|< 1$, which is true for sufficiently small $\Delta\tau$ and under the previous assumption that $\Delta\tau\ll\Delta t$, and hence may treat the summation as a geometric series.
	
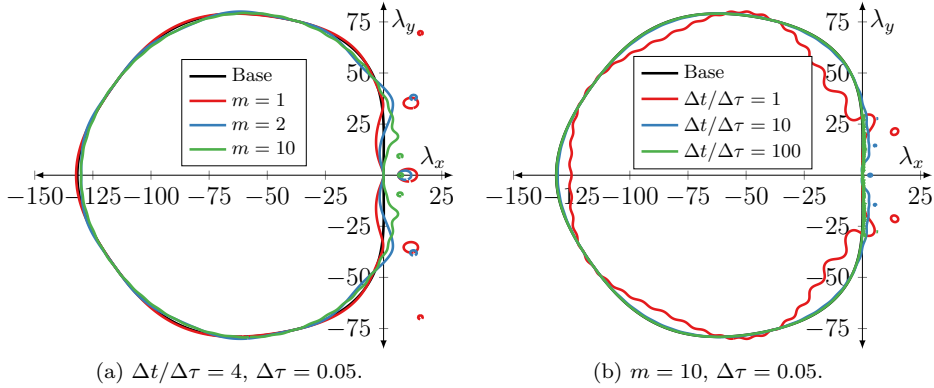
\begin{figure}[tbhp]
	    \centering
        \subfloat[$\Delta t/\Delta\tau=4$, $\Delta\tau = 0.05$.]{\label{fig:stability_region} \adjustbox{width=0.48\linewidth,valign=b}{	\begin{tikzpicture}
		\begin{axis}[name=plot1,xlabel={$\lambda_x$},ylabel={$\lambda_y$},
		    axis lines=middle,
            axis line style={latex-latex},
		    xtick={-150,-125,...,25},
		    xmin=-150,xmax=30,
		    ytick={-75,-50,...,75},
    		ymin=-85,ymax=85,
    		legend style={at={(0.5,0.68)},anchor=center,font=\small},
    		legend cell align={left},
    		style={font=\large}
    		]
    		
			\addplot[color=black, style={very thick}]
				table[x=x,y=y,col sep=comma,unbounded coords=jump]{./Figs/data/BDF_stability/OERK5dgp4.csv};
			\addlegendentry{Base}
			
			\addplot[color={Set1-A}, style={very thick}] table[x=x,y=y,col sep=comma,unbounded coords=jump]{./Figs/data/BDF_stability/OERK5dgp4_BDF2-1.csv};
			\addplot[color={Set1-A}, style={very thick},forget plot] table[x=x,y=y,col sep=comma,unbounded coords=jump]{./Figs/data/BDF_stability/OERK5dgp4_BDF2-1a.csv};
			\addplot[color={Set1-A}, style={very thick},forget plot] table[x=x,y=y,col sep=comma,unbounded coords=jump]{./Figs/data/BDF_stability/OERK5dgp4_BDF2-1b.csv};
			\addplot[color={Set1-A}, style={very thick},forget plot] table[x=x,y=y,col sep=comma,unbounded coords=jump]{./Figs/data/BDF_stability/OERK5dgp4_BDF2-1c.csv};
			\addplot[color={Set1-A}, style={very thick},forget plot]	table[x=x,y=y,col sep=comma,unbounded coords=jump]{./Figs/data/BDF_stability/OERK5dgp4_BDF2-1d.csv};
			\addplot[color={Set1-A}, style={very thick},forget plot]	table[x=x,y=y,col sep=comma,unbounded coords=jump]{./Figs/data/BDF_stability/OERK5dgp4_BDF2-1e.csv};
			\addlegendentry{$m=1$}
			
			\addplot[color={Set1-B}, style={very thick}] table[x=x,y=y,col sep=comma,unbounded coords=jump]{./Figs/data/BDF_stability/OERK5dgp4_BDF2-2.csv};
			\addplot[color={Set1-B}, style={very thick},forget plot] table[x=x,y=y,col sep=comma,unbounded coords=jump]{./Figs/data/BDF_stability/OERK5dgp4_BDF2-2a.csv};
			\addplot[color={Set1-B}, style={very thick},forget plot] table[x=x,y=y,col sep=comma,unbounded coords=jump]{./Figs/data/BDF_stability/OERK5dgp4_BDF2-2b.csv};
			\addplot[color={Set1-B}, style={very thick},forget plot] table[x=x,y=y,col sep=comma,unbounded coords=jump]{./Figs/data/BDF_stability/OERK5dgp4_BDF2-2c.csv};
			\addplot[color={Set1-B}, style={very thick},forget plot] table[x=x,y=y,col sep=comma,unbounded coords=jump]{./Figs/data/BDF_stability/OERK5dgp4_BDF2-2d.csv};
			\addplot[color={Set1-B}, style={very thick},forget plot] table[x=x,y=y,col sep=comma,unbounded coords=jump]{./Figs/data/BDF_stability/OERK5dgp4_BDF2-2e.csv};
			\addlegendentry{$m=2$}
			
			\addplot[color={Set1-C}, style={very thick}] table[x=x,y=y,col sep=comma,unbounded coords=jump]{./Figs/data/BDF_stability/OERK5dgp4_BDF2-10.csv};
			\addplot[color={Set1-C}, style={very thick},forget plot] table[x=x,y=y,col sep=comma,unbounded coords=jump]{./Figs/data/BDF_stability/OERK5dgp4_BDF2-10a.csv};
			\addplot[color={Set1-C}, style={very thick},forget plot] table[x=x,y=y,col sep=comma,unbounded coords=jump]{./Figs/data/BDF_stability/OERK5dgp4_BDF2-10b.csv};
			\addplot[color={Set1-C}, style={very thick},forget plot] table[x=x,y=y,col sep=comma,unbounded coords=jump]{./Figs/data/BDF_stability/OERK5dgp4_BDF2-10c.csv};
			\addlegendentry{$m=10$}
		\end{axis} 		
	\end{tikzpicture}}}
        \subfloat[$m=10$, $\Delta\tau=0.05$.]{\label{fig:dt_stability} \adjustbox{width=0.48\linewidth,valign=b}{	\begin{tikzpicture}
		\begin{axis}[name=plot1,xlabel={$\lambda_x$},ylabel={$\lambda_y$},
		    axis lines=middle,
            axis line style={latex-latex},
		    xtick={-150,-125,...,25},
		    xmin=-150,xmax=30,
		    ytick={-75,-50,...,75},
    		ymin=-85,ymax=85,
    		legend style={at={(0.5,0.68)},anchor=center,font=\small},
    		legend cell align={left},
    		style={font=\large}
    		]
    		
			\addplot[color=black, style={very thick}]
				table[x=x,y=y,col sep=comma,unbounded coords=jump]{./Figs/data/BDF_stability/OERK5dgp4.csv};
			\addlegendentry{Base}
			
			\addplot[color={Set1-A}, style={very thick}] table[x=x,y=y,col sep=comma,unbounded coords=jump]{./Figs/data/BDF_stability/OERK5dgp4_BDF2-dt-1.csv};
			\addplot[color={Set1-A}, style={very thick},forget plot] table[x=x,y=y,col sep=comma,unbounded coords=jump]{./Figs/data/BDF_stability/OERK5dgp4_BDF2-dt-1a.csv};
			\addplot[color={Set1-A}, style={very thick},forget plot] table[x=x,y=y,col sep=comma,unbounded coords=jump]{./Figs/data/BDF_stability/OERK5dgp4_BDF2-dt-1b.csv};
			\addplot[color={Set1-A}, style={very thick},forget plot] table[x=x,y=y,col sep=comma,unbounded coords=jump]{./Figs/data/BDF_stability/OERK5dgp4_BDF2-dt-1c.csv};
			\addlegendentry{$\Delta t/\Delta\tau = 1$}
			
			\addplot[color={Set1-B}, style={very thick}] table[x=x,y=y,col sep=comma,unbounded coords=jump]{./Figs/data/BDF_stability/OERK5dgp4_BDF2-dt-3.csv};
			\addplot[color={Set1-B}, style={very thick},forget plot] table[x=x,y=y,col sep=comma,unbounded coords=jump]{./Figs/data/BDF_stability/OERK5dgp4_BDF2-dt-3a.csv};
			\addplot[color={Set1-B}, style={very thick},forget plot] table[x=x,y=y,col sep=comma,unbounded coords=jump]{./Figs/data/BDF_stability/OERK5dgp4_BDF2-dt-3b.csv};
			\addplot[color={Set1-B}, style={very thick},forget plot] table[x=x,y=y,col sep=comma,unbounded coords=jump]{./Figs/data/BDF_stability/OERK5dgp4_BDF2-dt-3c.csv};
			\addplot[color={Set1-B}, style={very thick},forget plot] table[x=x,y=y,col sep=comma,unbounded coords=jump]{./Figs/data/BDF_stability/OERK5dgp4_BDF2-dt-3d.csv};
			\addplot[color={Set1-C}, style={very thick},forget plot] table[x=x,y=y,col sep=comma,unbounded coords=jump]{./Figs/data/BDF_stability/OERK5dgp4_BDF2-dt-3e.csv};
			\addlegendentry{$\Delta t/\Delta\tau = 10$}
			
			\addplot[color={Set1-C}, style={very thick}] table[x=x,y=y,col sep=comma,unbounded coords=jump]{./Figs/data/BDF_stability/OERK5dgp4_BDF2-dt-4.csv};
			\addplot[color={Set1-C}, style={very thick},forget plot] table[x=x,y=y,col sep=comma,unbounded coords=jump]{./Figs/data/BDF_stability/OERK5dgp4_BDF2-dt-4a.csv};
			\addplot[color={Set1-C}, style={very thick},forget plot] table[x=x,y=y,col sep=comma,unbounded coords=jump]{./Figs/data/BDF_stability/OERK5dgp4_BDF2-dt-4b.csv};
			\addplot[color={Set1-C}, style={very thick},forget plot] table[x=x,y=y,col sep=comma,unbounded coords=jump]{./Figs/data/BDF_stability/OERK5dgp4_BDF2-dt-4c.csv};
			\addplot[color={Set1-C}, style={very thick},forget plot] table[x=x,y=y,col sep=comma,unbounded coords=jump]{./Figs/data/BDF_stability/OERK5dgp4_BDF2-dt-4d.csv};
			\addplot[color={Set1-C}, style={very thick},forget plot] table[x=x,y=y,col sep=comma,unbounded coords=jump]{./Figs/data/BDF_stability/OERK5dgp4_BDF2-dt-4e.csv};
			\addplot[color={Set1-C}, style={very thick},forget plot] table[x=x,y=y,col sep=comma,unbounded coords=jump]{./Figs/data/BDF_stability/OERK5dgp4_BDF2-dt-4f.csv};
			\addplot[color={Set1-C}, style={very thick},forget plot] table[x=x,y=y,col sep=comma,unbounded coords=jump]{./Figs/data/BDF_stability/OERK5dgp4_BDF2-dt-4g.csv};
			\addplot[color={Set1-C}, style={very thick},forget plot] table[x=x,y=y,col sep=comma,unbounded coords=jump]{./Figs/data/BDF_stability/OERK5dgp4_BDF2-dt-4h.csv};
			\addplot[color={Set1-C}, style={very thick},forget plot] table[x=x,y=y,col sep=comma,unbounded coords=jump]{./Figs/data/BDF_stability/OERK5dgp4_BDF2-dt-4i.csv};
			\addplot[color={Set1-C}, style={very thick},forget plot] table[x=x,y=y,col sep=comma,unbounded coords=jump]{./Figs/data/BDF_stability/OERK5dgp4_BDF2-dt-4j.csv};
			\addplot[color={Set1-C}, style={very thick},forget plot] table[x=x,y=y,col sep=comma,unbounded coords=jump]{./Figs/data/BDF_stability/OERK5dgp4_BDF2-dt-4k.csv};
			\addplot[color={Set1-C}, style={very thick},forget plot] table[x=x,y=y,col sep=comma,unbounded coords=jump]{./Figs/data/BDF_stability/OERK5dgp4_BDF2-dt-4l.csv};
			\addplot[color={Set1-C}, style={very thick},forget plot] table[x=x,y=y,col sep=comma,unbounded coords=jump]{./Figs/data/BDF_stability/OERK5dgp4_BDF2-dt-4m.csv};
			\addplot[color={Set1-C}, style={very thick},forget plot] table[x=x,y=y,col sep=comma,unbounded coords=jump]{./Figs/data/BDF_stability/OERK5dgp4_BDF2-dt-4n.csv};
			\addplot[color={Set1-C}, style={very thick},forget plot] table[x=x,y=y,col sep=comma,unbounded coords=jump]{./Figs/data/BDF_stability/OERK5dgp4_BDF2-dt-4o.csv};
			\addplot[color={Set1-C}, style={very thick},forget plot] table[x=x,y=y,col sep=comma,unbounded coords=jump]{./Figs/data/BDF_stability/OERK5dgp4_BDF2-dt-4p.csv};
			\addplot[color={Set1-C}, style={very thick},forget plot] table[x=x,y=y,col sep=comma,unbounded coords=jump]{./Figs/data/BDF_stability/OERK5dgp4_BDF2-dt-4q.csv};
			\addplot[color={Set1-C}, style={very thick},forget plot] table[x=x,y=y,col sep=comma,unbounded coords=jump]{./Figs/data/BDF_stability/OERK5dgp4_BDF2-dt-4r.csv};
			\addplot[color={Set1-C}, style={very thick},forget plot] table[x=x,y=y,col sep=comma,unbounded coords=jump]{./Figs/data/BDF_stability/OERK5dgp4_BDF2-dt-4s.csv};
			\addplot[color={Set1-C}, style={very thick},forget plot] table[x=x,y=y,col sep=comma,unbounded coords=jump]{./Figs/data/BDF_stability/OERK5dgp4_BDF2-dt-4t.csv};
			\addplot[color={Set1-C}, style={very thick},forget plot] table[x=x,y=y,col sep=comma,unbounded coords=jump]{./Figs/data/BDF_stability/OERK5dgp4_BDF2-dt-4u.csv};
			\addplot[color={Set1-C}, style={very thick},forget plot] table[x=x,y=y,col sep=comma,unbounded coords=jump]{./Figs/data/BDF_stability/OERK5dgp4_BDF2-dt-4v.csv};
			\addplot[color={Set1-C}, style={very thick},forget plot] table[x=x,y=y,col sep=comma,unbounded coords=jump]{./Figs/data/BDF_stability/OERK5dgp4_BDF2-dt-4w.csv};
			\addplot[color={Set1-C}, style={very thick},forget plot] table[x=x,y=y,col sep=comma,unbounded coords=jump]{./Figs/data/BDF_stability/OERK5dgp4_BDF2-dt-4x.csv};
			\addlegendentry{$\Delta t/\Delta\tau = 100$}

			
		\end{axis} 		
	\end{tikzpicture}}}
        \caption{\label{fig:ERK5dgp4_BDF2_stab}Stability regions for OERK5-DGp4 and BDF2.}
    \end{figure}

	The contours of unity amplification factor representing the stability limit are shown in \cref{fig:ERK5dgp4_BDF2_stab} for BDF2 coupled to an ERK scheme. The ERK scheme applied was an optimised 5 stage scheme from the work of Vermeire~\etal~\cite{Vermeire2019}, where the ERK stability region was tuned to match the set of eigenvalues produced by $p=4$ nodal DG spatial scheme for advection. This scheme will be denoted as \emph{OERK5-DGp4}. It was posited that the schemes would provide the optimal stability region when using this spatial scheme with dual-time stepping for implicit calculations. To produce a stability region, it was  necessary to set $\Delta t$ and $\Delta\tau$, which for these contours take the value of $0.2$ and $0.05$, respectively. Then, eigenvalues can be applied to the system to find the unity contour using $\lambda=\lambda_x+\imath\lambda_y$.
	
     As is demonstrated here, the coupling of the implicit method to the pseudo-time integrator causes the stability region to change with the number of iterations, with both local contractions and expansions observed. The stability region of the ERK scheme without coupling to an implicit method is also shown in \cref{fig:ERK5dgp4_BDF2_stab} for reference. Additionally, the stability is further deformed by variations to the ratio $\Delta t/\Delta \tau$, \cref{fig:dt_stability}, therefore complicating the design of optimal ERK schemes. Further investigation of the stability of the dual-time system was performed by Chiew~\etal~\cite{Chiew2016} where a more exhaustive study of implicit schemes is given. 
	
\subsection{\emph{p}-multigrid}
    To accelerate the convergence of the solution towards a \newline pseudo-time steady state, the \emph{p}-multigrid methodology has proven to be effective for spectral element methods such as FR~\cite{Loppi2018}. The aim of the method is to restrict the solution to coarse grid levels, apply smoothing there, and subsequently propagate corrections from the coarser levels to the finer levels. We will now outline the techniques of \emph{p}-multigrid applied to the system already described. From the work of the previous section, the residual after $M$ pseudo time steps is
	\begin{equation}\label{eq:ps_residual}
		T_{p,M}u_{n+1,0} = -\lambda\bigg[P^M + \frac{1-P^M}{1-P}C\bigg]u_{n+1,0} - \frac{1}{\Delta t B_0}\bigg[u_{n+1,M} - \sum^{s-1}_{l=0}B_{1+l}u_{n-l}\bigg].
	\end{equation}
	For the finest stage, of degree $p$, the deficit is defined as
	\begin{equation}
		d_p = -T_{p,M}.
	\end{equation}
	The deficit and residual source terms for the lower order stages are subsequently defined as
	\begin{subequations}
		\begin{align}
			u_{i-1,n+1,0} =&\: \rho_{i-1}(u_{i,n+1,M}), \\
			d_{i-1} =&\: \rho_{i-1}(d_{i}), \\
			r_{i-1} =&\: T_{i-1,M}u_{i-1,n+1,0} + d_{i-1},
		\end{align}
	\end{subequations}
	where $r_{i-1}$ is the deficit residual source term that is applied in the calculation of $u_{i-1,n+1,M}$, to be shown momentarily. The restriction operator, $\rho_i(.)$, is taken to be the same for the solution and deficit and is defined as 
	\begin{equation}
	    \langle\rho_k(u)-u,\phi_i\rangle_{L_2} = 0,
	\end{equation}
	for some polynomial basis $\phi_i$, which we will take to be the orthogonal Legendre basis. For the linear case to be considered here, this choice does not restrict the generality of the results; however, otherwise this choice is justified by being a polynomial basis for $L_2$ polynomial projection with unit measure. When defined within a nodal or collocation spatial method, the inner product will require approximation for which we use quadrature rules such as Gauss--Legendre. 
	
	The prolongation and correction of the $i^\mathrm{th}$ level based on the $i-1^\mathrm{th}$ is then
	\begin{subequations}
		\begin{align}
			\Delta_{i} =&\: v_{i,n+1,0} - v_{i,n+1,M},\\
			\Delta_{i+1} =&\: \pi_{i+1}(\Delta_i), \\
			v_{i+1,n+1,0} =&\: u_{i+1,n+1,M} + c_{i+1},
		\end{align}
	\end{subequations}
	where $v$ is used to indicate the new solution on the prolongation steps. If at a local minima in a \emph{p}-multigrid cycle, $v_{i,n+1,0}$ is taken to be $u_{i,n+1,0}$. Furthermore, the prolongation operator $\pi_i$ is defined such that given $u_k\in\mathbb{P}^k$ and $x_{k,i}\in\{x_0,\dots,x_k\}$, then
	\begin{equation}
	    \pi_{k+1}(u_k) \in \mathbb{P}^{k+1} \quad \mathrm{given} \quad \pi_{k+1}(u_k)(x_{k+1,i}) = u_k(x_{k+1,i}).
	\end{equation}
	We now wish to incorporate the multi-grid residual source term $r_q$ such that the modified pseudo-time update equation may be defined, which manifests straightforwardly in the ERK steps as
	\begin{equation*}
		q_i = \lambda \Big(u_{n+1,m} + \Delta\tau\sum^{i-1}_{j=1}a_{i,j}k_j\Big) - \frac{1}{\Delta tB_0}\Big(u_{n+1,m} + \sum^{s-1}_{l=0}B_{l+1}u_{n-l}\Big) - r_q,
	\end{equation*}	 
	and hence
	\begin{equation}\label{eq:ps_rsource}
		u_{n+1,m+1} = Pu_{n+1,m}-C\sum^{s-1}_{l=0}B_{l+1}u_{n-l}  - \Big(\underbrace{\Delta\tau\mathbf{b}^T(\mathbf{I} - \lambda\Delta\tau\mathbf{A})^{-1}\mathbf{e}}_{K}\Big)r_q.
	\end{equation}

\section{Fourier Analysis}\label{sec:fourier}
    We have so far presented the techniques to construct implicit temporal integration applied to ODEs, \cref{eq:linear_pseduo} and demonstrated the effect of pseudo-stepping on time integration stability. We now wish to use the flux reconstruction scheme for spatial differentiation to provide the eigenvalues. With this complete system, not only can the coupled stability be studied, but it provides a means to calculate the analytic error which may  inform cycle construction. In order to generalise the analysis, we will consider the Fourier analysis of the linear advection-diffusion equation with a modified Bloch trial solution
    \begin{subequations}\label{eq:lad}
        \begin{align}
            \px{u}{t} + \px{u}{x} &= \mu\pxi{2}{u}{x}, \\
            u &= \exp{\big(\imath\left(kx-\omega t\right)\big)},\\
		    \mathbf{u}_n &= \exp{\big(\imath (\mathbf{x}-\omega n\Delta t)\big)},
		\end{align}
	\end{subequations}
    where $k$ is the wavenumber, $\omega=k(1-\imath \mu k)$ is the angular frequency, and $\imath = \sqrt{-1}$. We will now construct the spatial derivatives via the FR methodology~\cite{Huynh2007,Vincent2010} in one-dimension, which in the linear case---with the Bloch wave solution---may be defined as 
	\begin{subequations}\label{eq:Q_def}
	    \begin{align}
		    \px{\mathbf{u}_i}{x} = \mathbf{Q}_a\mathbf{u}_i =&\; \frac{2}{h}\Big(\exp{(-\imath kh)}\mathbf{C}_{-} + \mathbf{C}_0 + \exp{(\imath kh)}\mathbf{C}_+\Big)\mathbf{u}_i,\\
		    \pxi{2}{\mathbf{u}_i}{x} = \mathbf{Q}_d\mathbf{u}_i =&\; 
		    \frac{4}{h^2}\Big(\exp{(-\imath 2kh)}\mathbf{B}_{-2} + 
		    \exp{(-\imath kh)}\mathbf{B}_{-} + \mathbf{B}_0 + \\
		    &\quad\quad\:\exp{(\imath kh)}\mathbf{B}_{+} + 
		   \exp{(2\imath kh)}\mathbf{B}_{+2}\Big)\mathbf{u}_i, \notag
		\end{align}
	\end{subequations}
	for linearly transformed elements on a uniform grid with spacing $h$. Further details on the  operator definitions can be found in \cref{app:fr_op}.  During the FR method, a common interface flux and a common interface value is calculated. We will use $\alpha=(\alpha_a,\alpha_d)$ to denote the degree of upwinding in the advection and diffusion calculations, with $\alpha=1$ being fully upwinded and $\alpha=0.5$ being centrally differenced. 
	
	If the FR scheme represented by $\mathbf{Q}$ is full rank, i.e., none of the solution points are collocated and $k\ne 0$, then $\mathbf{Q}$ can be diagonalised as
	\begin{equation}\label{eq:q_diag}
		\mathbf{Q} = -\mathbf{Q}_a + \mu\mathbf{Q}_d = \imath k\mathbf{W\Lambda}_Q\mathbf{W}^{-1}. 
	\end{equation}
	which demonstrates that FR has the capacity for a solution comprised of $p+1$ unique eigenvalues. To now apply FR as the source of the eigenvalues to the integration scheme, we first use a result of Ketcheson~\etal~\cite{Ketcheson2012}, where it is possible to write the stability polynomial of a temporal integration method with $r$ steps as
	\begin{equation}
		P\bigg(\lambda\Delta\tau,\frac{\Delta\tau}{\Delta t}\bigg) = \sum^{r}_{j=0}\gamma_j\bigg(\frac{\Delta\tau}{\Delta t}\bigg)(\lambda\Delta\tau)^j.
	\end{equation}
	for an $r$ stage RK scheme. Hence, we can define the partial pseudo update equation as
	\begin{equation}
		\mathbf{P} = \sum^{r}_{j=0}\gamma_j(\Delta\tau\mathbf{Q})^j = \mathbf{W}\Bigg[\sum^{r}_{j=0}\gamma_j(\imath k\Delta\tau\mathbf{\Lambda}_Q\big)^j\Bigg]\mathbf{W}^{-1}.
	\end{equation}
	The BDF source term is also a function of $\lambda\Delta\tau$ and can similarly be found in terms of $\mathbf{Q}$ using a polynomial fit of $C$ as in \cref{eq:ps1}. Hence,
	\begin{equation}
		C = \sum^{r-1}_{j=0}\kappa_j(\lambda\Delta\tau)^j, \quad \mathrm{and} \quad \mathbf{C} = \mathbf{W}\Bigg[\sum^{r-1}_{j=0}\kappa_j\big(\imath k\Delta\tau\mathbf{\Lambda}_Q\big)^j\Bigg]\mathbf{W}^{-1}.
	\end{equation}
	The full pseudo-time update equation is then
	\begin{equation}
		\mathbf{u}_{n+1,m+1} = \mathbf{P}\mathbf{u}_{n+1,m} - \mathbf{C}\sum^{s-1}_{l=0}B_{l+1}\mathbf{u}_{n-l}.
	\end{equation}
	To confirm $\mathbf{P}$ and $\mathbf{C}$ are correctly defined the following relation should hold
	\begin{equation}\label{eq:sum_update}
		\mathbf{P} + \mathbf{C} = \mathbf{R},
	\end{equation}
	for the ERK update matrix, $\mathbf{R}$. With $\mathbf{P}$ and $\mathbf{C}$ defined, the $M^\mathrm{th}$ value can be expressed in terms of the initial value of the pseudo-stepping as
	\begin{equation}\label{eq:update}
		    \mathbf{u}_{n+1,M} = \mathbf{P}^M\mathbf{u}_{n+1,0} - (\mathbf{I} -\mathbf{P})^{-1}(\mathbf{I} - \mathbf{P}^M)
		        \mathbf{C}\Bigg(\sum^{s-1}_{l=0}B_{l+1}\mathbf{u}_{n-l}\Bigg).
	\end{equation}
	Again simplification was made through a geometric series and its matrix analogue, which has the generalised assumption that the spectral radius of $\mathbf{P}$ is less than unity, i.e., $\rho(\mathbf{P}) \leqslant 1$. This can be verified for suitable pseudo-time steps coupled to the previous assumption that $\Delta\tau \ll \Delta t$. 
	
	The dual-time update may then be written as
	\begin{subequations}
		\begin{align}
			\mathbf{u}_{n+1,M} =&\: \left[\mathbf{P}^M - (\mathbf{I} -\mathbf{P})^{-1}(\mathbf{I} - \mathbf{P}^M)\mathbf{C}\left(\sum^{s-1}_{l=0}B_{l+1}\exp{\left(\imath \omega l\Delta t\right)} \right)\right]\mathbf{u}_{n+1,0}, \\
			\mathbf{u}_{n+1,M} =&\: \mathbf{R}_M\mathbf{u}_{n+1,0}.
		\end{align}	
	\end{subequations}	
	The previous solution needed for the BDF source term is taken as the analytic solution from \cref{eq:lad}, which is consistent with a time history of fully converge solutions. Due to imposing a discretiation on the solution the system has a Nyquist limit on the maximum wavenumber, which due to the coupled space-time is
	\begin{equation}\label{eq:nyquist}
	    k_\mathrm{Nq} = \min{\left(\frac{\pi}{\Delta t},\frac{(p+1)\pi}{h}\right)} \quad \mathrm{and}\quad \hat{k}=\frac{\pi k}{k_\mathrm{Nq}}.
	\end{equation}
	with $\hat{k}$ being the normalised wavenumber.
	
	The exact solution from the applied Bloch wave can be projected into the solution space of FR using the eigenvectors of $\mathbf{Q}$ to obtain the vector of mode weights, $\pmb{\beta}$, via
	\begin{equation}
	    \mathbf{u}_0 = \exp{\left(\imath kx_j\right)}\mathbf{W}\pmb{\beta}.
	\end{equation}
	This may then be substituted into \cref{eq:update} to give the fully discrete error, written as
	\begin{subequations}
		\begin{align}
			\mathbf{e} =&\: \mathbf{u}_{n+1,M} - \mathbf{u}_{n+1}, \\
			=&\: \exp{\big(\imath (kx_j-\omega n\Delta t)\big)}\left(\mathbf{R}_M -\exp{(-\imath \omega\Delta t)}\mathbf{I}\right)\mathbf{W}\pmb{\beta}.
		\end{align}	
	\end{subequations}
	
	\begin{figure}[tbhp]
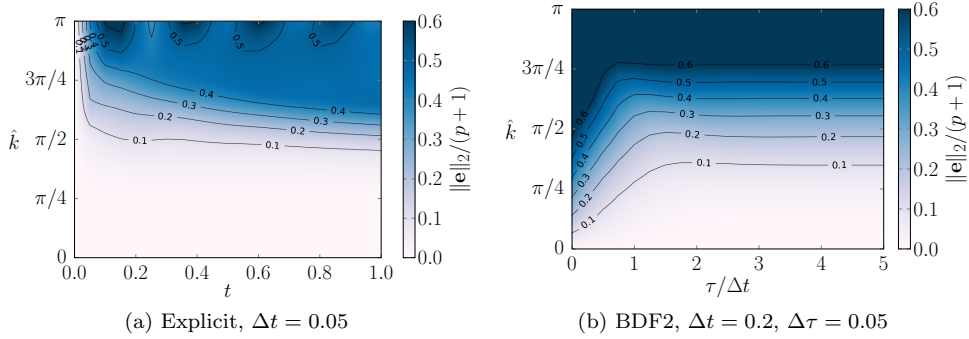

	    \centering
        \subfloat[Explicit, $\Delta t=0.05$]{\label{fig:frdg4_error} \adjustbox{width=0.48\linewidth,valign=b}{\input{./Figs/FRDG4_explicit_error.tex}}}
        ~
        \subfloat[BDF2, $\Delta t = 0.2$, $\Delta\tau=0.05$]{\label{fig:frdg4_error_bdf} \adjustbox{width=0.48\linewidth,valign=b}{\input{./Figs/FRDG4_BDF2_error_tt4.tex}}}
        \caption{\label{fig:frdg4_error_comp}Error comparison for FRDG, $p=4$, and SSPRK3 explicit scheme, with and without dual-time stepping for pure advection.}
    \end{figure}
    
    The evolution of the Euclidean norm of the error calculated using this method is shown in \cref{fig:frdg4_error_comp}, where comparison is made between the use of an explicit scheme and dual-time stepping for FR with upwinded interfaces. The physical time step size in the dual-time error was chosen such that the temporal and spatial Nyquist wavenumbers were equivalent. It is evident that at low wavenumbers the error is equivalent, but at high wavenumbers, the dispersion and dissipation associated with the scheme causes a modification to the pseudo-time steady state, and so large errors are observed.
	
	If we now look to characterise the maximum time step sizes for the explicit and coupled system, due to the presence of source terms in the update equation, the traditional von Neumann stability criteria has to be modified. Therefore, the set of stable values of $\Delta\tau$ may be defined as 
	\begin{equation}
	    \Delta\tau_\mathrm{stable}(\Delta t) = \Big\{\Delta\tau\in\mathbb{R}_+ : \: \rho\big(\mathbf{R}_M(\Delta\tau,\Delta t)\big) \leqslant \Big|\sum^{s-1}_{l=0}B_{l+1}\exp{\big(\imath \omega l\Delta t\big)}\Big|\Big\}.
	\end{equation}
	Hence, the maximum stable step size is $\Delta\tau_\mathrm{max} = \sup\Delta\tau_\mathrm{stable}$. We will also define  $\Delta\tau_\mathrm{max,A}$ to signify the maximum step size for pure advection.
	
	 \begin{figure}[tbhp]
	    \centering
	    \subfloat[Maximum time step size for advection-diffusion FRDG with explicit SSPRK3 temporal integration, and $\alpha=(1,0.5)$.]
	    {\label{fig:FRDG_ad_cfl}\adjustbox{width=0.48\linewidth,valign=b}{\begin{tikzpicture}
		\begin{axis}[name=plot1,xlabel={$\mu$},ylabel={$\Delta t_\mathrm{max}$},
            axis line style={latex-latex},
            axis y line=middle,
            axis x line=left,
            xmode=log,
		    xtick={1e-3,1e-2,1e-1,1e0,1e1},
		    xmin=1e-3,xmax=1e1,
            ymode=log,
		    ytick={1e-3,1e-2,1e-1,1e0},
    		ymin=1e-3,ymax=1,
    		legend style={at={(0.05,0.05)},anchor=south west,font=\small},
    		style={font=\large}
    		]

            \addplot[color={Set1-A}, style={very thick}]
				table[x=mu,y=cfl1,col sep=comma,unbounded coords=jump]{./Figs/data/FR_ssprk3_up_cent_ad_cfl.csv};
			\addlegendentry{$p=1$}
			
            \addplot[color={Set1-B}, style={very thick}]
				table[x=mu,y=cfl2,col sep=comma,unbounded coords=jump]{./Figs/data/FR_ssprk3_up_cent_ad_cfl.csv};
			\addlegendentry{$p=2$}
			
            \addplot[color={Set1-C}, style={very thick}]
				table[x=mu,y=cfl3,col sep=comma,unbounded coords=jump]{./Figs/data/FR_ssprk3_up_cent_ad_cfl.csv};
			\addlegendentry{$p=3$}
			
            \addplot[color={Set1-D}, style={very thick}]
				table[x=mu,y=cfl4,col sep=comma,unbounded coords=jump]{./Figs/data/FR_ssprk3_up_cent_ad_cfl.csv};
			\addlegendentry{$p=4$}
			
            \addplot[color={Set1-E}, style={very thick}]
				table[x=mu,y=cfl5,col sep=comma,unbounded coords=jump]{./Figs/data/FR_ssprk3_up_cent_ad_cfl.csv};
			\addlegendentry{$p=5$}

    \end{axis}    		
\end{tikzpicture}}}
        ~
        \subfloat[{Maximum pseudo time step size for upwinded advection with FRDG $p=4$, BDF2, and SSPRK3 with pseudo step number $m$ and $k\in(0,k_\mathrm{Nq}]$.}]
        {\label{fig:FRDG_BDF_CFL}\adjustbox{width=0.48\linewidth,valign=b}{\begin{tikzpicture}
		\begin{axis}[name=plot1,xlabel={$\Delta t$},ylabel={$\Delta \tau_\mathrm{max}$},
            axis line style={latex-latex},
            axis y line=middle,
            axis x line=left,
            xmode=log,
		    xtick={1e-3,1e-2,1e-1,1e0,1e1},
		    xmin=1e-3,xmax=1e1,
            ytick={0.02,0.04,...,0.1},
		    ylabel style={at={(0.75,0.95)},anchor=east},
    		ymin=0,ymax=0.12,
    		yticklabel pos=right,
    		legend style={at={(0.02,0.72)},anchor=north west,font=\small},
    		legend cell align={left},
    		style={font=\large}
    		]

            \addplot[color={black}, style={dashed}] coordinates {
		        (1e-3,0.08968)
		        (1e1,0.08968)};
			\addlegendentry{Explicit}

            \addplot[color={Set1-A}, style={very thick}]
				table[x=dt,y=dtaumax,col sep=comma,unbounded coords=jump]{./Figs/data/FRp4DG_SSPRK3_1BDF2_CFL.csv};
			\addlegendentry{$m=1$}
			
			\addplot[color={Set1-B}, style={very thick}]
				table[x=dt,y=dtaumax,col sep=comma,unbounded coords=jump]{./Figs/data/FRp4DG_SSPRK3_2BDF2_CFL.csv};
			\addlegendentry{$m=2$}
			
			\addplot[color={Set1-C}, style={very thick}]
				table[x=dt,y=dtaumax,col sep=comma,unbounded coords=jump]{./Figs/data/FRp4DG_SSPRK3_10BDF2_CFL.csv};
			\addlegendentry{$m=10$}

    \end{axis}    		
\end{tikzpicture}}}
        \caption{\label{fig:CFL}Maximum time step size for some configurations for FRDG.}
    \end{figure}
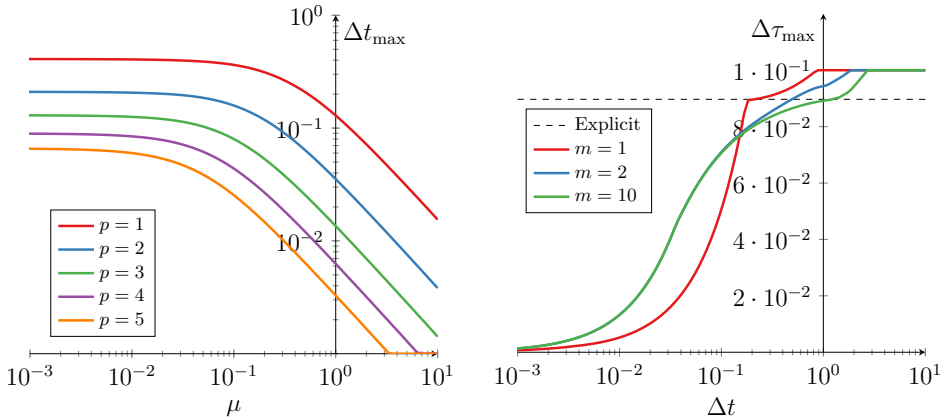
    
    \cref{fig:FRDG_ad_cfl} presents the CFL limits of the explicit system, i.e. without dual-time stepping, and makes it clear that for all orders the absolute value of the maximal explicit time-step becomes severely limited for low Reynolds numbers. Turning to the coupled system, the effect of the physical time size on the maximum pseudo-step size is presented in \cref{fig:FRDG_BDF_CFL}. Interestingly, it can be seen that the first pseudo-step has a more restrictive maximum step size, and from the error in \cref{fig:frdg4_error_bdf}, this can be attributed to the contraction being highest initially. Therefore, to prevent instabilities initially entering the solution, smaller pseudo-time steps are required at first. From \cref{fig:ERK5dgp4_BDF2_stab}, as the ratio $\Delta t/\Delta\tau$ is reduced the stability region is reduced and this is seen here in the CFL limit. For $\Delta t>0.2$, the physical-time dominates the Nyqusit limit, and it is around this point at which a sharp change in the $m=1$ case is seen. After this point, as $\Delta t$ is continually increased, the range of wavenumber decreases, and the stability is observed to increase. This is concurrent with the initial error in the BDF approximation being largest at highest wavenumbers, with further iterations this behaviour is not seen as the poor initial approximation of the temporal derivative from BDF---due to the use of $u_{n+1,0}=u_n$---is quickly rectified. 
    
    \subsection{\emph{p}-Multigrid}
		A key component of the multigrid methodology is the residual which was defined in \cref{eq:ps_residual} for the $M^\mathrm{th}$ pseudo step in terms of the zeroth step. Through applying the FR operator for the spatial discretisation, we may write the residual of $u_{n+1,M}$ as:
		\begin{subequations}\label{eq:pseudo_res}
			\begin{align}
				\px{\mathbf{u}_{i,n+1,M}}{\tau} =&\: -\bigg(\frac{\mathbf{I}}{\Delta tB_0} + \mathbf{Q}_i\bigg)\mathbf{u}_{i,n+1,M} - \underbrace{\frac{1}{\Delta t B_0}\sum_{l=0}^{s-1}B_{l+1}\exp{\big(\imath \omega l\Delta t\big)}}_{C_B}\mathbf{u}_{i,n+1,0}  \\
				=&\: \mathbf{T}_{i}\mathbf{u}_{i,n+1,M}  - C_B\mathbf{u}_{i,n+1,0}.
			\end{align}
		\end{subequations}	
		For FR \emph{p}-multigrid, the restriction and prolongation matrices can be straightforwardly defined modally as
		\begin{equation}
		   	\hat{\pmb{\rho}}_i = \hat{\mathbf{I}} \quad \mathrm{and} \quad \hat{\pmb{\pi}}_i = \hat{\mathbf{I}}^T,
		\end{equation}
		which can be projected to a nodal representation by using the Vandermonde matrix and the appropriate solution points for the degree.
			
		It should be noted again that $\mathbf{r}_p = 0$. To proceed, \cref{eq:ps_rsource} has to converted to a matrix representation and so the procedure applied to $C$ is applied to $K$ using $\mathbf{Q}$. This leads to the update equation
		\begin{equation}
			\mathbf{u}_{i,n+1,m+1} = \mathbf{P}_i\mathbf{u}_{i,n+1,m} - \mathbf{C}_i\sum^{s-1}_{l=0}B_{l+1}\mathbf{u}_{i,n-l} - \mathbf{K}_i\mathbf{r}_i,
		\end{equation}
		which similarly may be defined at the $M^\mathrm{th}$ step as
		\begin{subequations}
			\begin{align}
				\mathbf{u}_{i,n+1,M} =&\: \mathbf{P}^M_i\mathbf{u}_{i,n+1,0} - \bigg[\sum^{M-1}_{m=0}\mathbf{P}^m_i\bigg]\bigg(\mathbf{C}_i\sum^{s-1}_{l=0}B_{l+1}\mathbf{u}_{i,n-l} + \mathbf{K}_i\mathbf{r}_i\bigg),\\
				\mathbf{u}_{i,n+1,M} =&\: \mathcal{S}(M,\mathbf{P},\mathbf{C},\mathbf{K},\mathbf{r}_i,\mathbf{u}_{i,n+1,0},\mathbf{u}_i).
			\end{align}
		\end{subequations}
		
		\begin{figure}[tbhp]
	        \centering
            \subfloat[][One level \emph{V}-cycle.]{\label{fig:simple_v} \adjustbox{width=0.3\linewidth,valign=b}{\begin{tikzpicture}
	\begin{axis}[%
    	width=0.2\textwidth,
    	height=0.2\textwidth,
        scale only axis,
        xmin=-0.5, xmax=2.5,
        xticklabels={,,},
        ymin=-1, ymax=0,
        ytick={-1,0},
        yticklabels={$p-1$,$p$},
        ylabel style={yshift=0.4em},
        ymajorgrids,
        line width=2.0pt,
        mark size=4.0pt,
        axis line style={draw=none},
        tick style={draw=none}
        ]
        \addplot[color=black,thin,solid,mark=*,mark options={solid},mark size=1.3pt] coordinates{(0,0) (1,-1) (2,0)};
	\end{axis}
\end{tikzpicture}}}
            ~
            \subfloat[][Multi-level $\mathrm{\emph{W}}_{p-2}\,$-cycle. ]{\label{fig:simple_w}\adjustbox{width=0.3\linewidth,valign=b}{\begin{tikzpicture}
	\begin{axis}[%
    	width=0.2\textwidth,
    	height=0.2\textwidth,
        scale only axis,
        xmin=-2, xmax=14,
        xticklabels={,,},
        ymin=-4, ymax=0,
        ytick={-4,-3,-2,-1,0},
        yticklabels={$p-4$,$p-3$,$p-2$,$p-1$,$p$},
        ylabel style={yshift=0.4em},
        ymajorgrids,
        line width=2.0pt,
        mark size=4.0pt,
        axis line style={draw=none},
        tick style={draw=none}
        ]
        \addplot[color=black,thin,solid,mark=*,mark options={solid},mark size=1.3pt] coordinates{(0,0) 
                                                                            (1,-1)
                                                                            (2,-2)
                                                                            (3,-3)
                                                                            (4,-4)
                                                                            (5,-3)
                                                                            (6,-2)
                                                                            (7,-3)
                                                                            (8,-4)
                                                                            (9,-3)
                                                                            (10,-2)
                                                                            (11,-1)
                                                                            (12,0)};
	\end{axis}
\end{tikzpicture}}}
            ~
            \subfloat[][Multi-level $\mathrm{\emph{V}}_\mathrm{AP}$-cycle. ]{\label{fig:ap_v}\adjustbox{width=0.3\linewidth,valign=b}{\begin{tikzpicture}
	\begin{axis}[%
    	width=0.2\textwidth,
    	height=0.2\textwidth,
        scale only axis,
        xmin=-1, xmax=13,
        xticklabels={,,},
        ymin=-4, ymax=0,
        ytick={-4,-3,-2,-1,0},
        yticklabels={$p-4$,$p-3$,$p-2$,$p-1$,$p$},
        ylabel style={yshift=0.4em},
        ymajorgrids,
        line width=2.0pt,
        mark size=4.0pt,
        axis line style={draw=none},
        tick style={draw=none},
    	style={font=\small},
        ]
        \addplot[color=black,thin,solid,mark=*,mark options={solid},mark size=1.3pt] coordinates{(0,0) 
                                                                                                 (1,-1)
                                                                                                 (2,-2)
                                                                                                 (3,-3)
                                                                                                 (4,-4)
                                                                                                 (5,-4)
                                                                                                 (6,-3)
                                                                                                 (7,-3)
                                                                                                 (8,-2)
                                                                                                 (9,-2)
                                                                                                 (10,-1)
                                                                                                 (11,-1)
                                                                                                 (12,0)};
	\end{axis}
\end{tikzpicture}}}
            \caption{\label{fig:cycle_config}\emph{p}-Multigrid cycle configuration.}
        \end{figure}
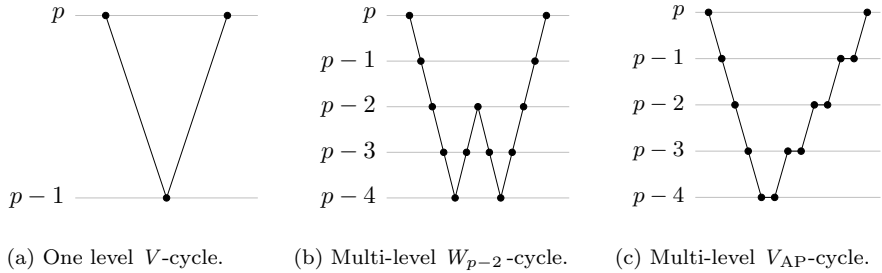
		
		We will now begin by defining the steps in a simple \emph{p}-multigrid \emph{V}-cycle. Diagrammatically, this is shown in \cref{fig:simple_v} with the steps presented in \cref{tab:simplev_steps}. This procedure may be generalised to arbitrary cycles such as \cref{tab:genv_steps} for \emph{V}-cycles with multiple stages.
		
        \begingroup
		\setlength{\tabcolsep}{2pt}
		\begin{table}[tbhp]
		    \overfullrule=0pt
		    \centering
			\caption{\label{tab:simplev_steps}Simple $V$-cycle steps.}
			\begin{tabular}{rclr}
					\toprule
					$\mathbf{u}_{p,n+1,M}$ &$=$& $\mathbf{R}_{p,M}\mathbf{u}_{p,n+1,0}$ & \rdelim\}{5}{0mm}[{\rotatebox[origin=c]{-90}{Restriction}}]\\
					$\mathbf{d}_{p}$ &$=$& $\mathbf{r}_p - (\mathbf{T}_{p}\mathbf{u}_{p,n+1,M} + C_B\mathbf{u}_{p,n+1,0})$ \\
					$\mathbf{u}_{p-1,n+1,0}$ &$=$& $\pmb{\rho}_{p-1}\mathbf{u}_{p,n+1,M}$\\
					$\mathbf{d}_{p-1}$ &$=$& $\pmb{\rho}_{p-1}\mathbf{d}_p$ \\
					$\mathbf{r}_{p-1}$ &$=$& $\mathbf{T}_{p-1}\mathbf{u}_{p-1,n+1,0} + C_B\pmb{\rho}_{p-1}\mathbf{u}_{p,n+1,0} + \mathbf{d}_{p-1}$ \\
					\midrule
					$\mathbf{u}_{p-1,n+1,M}$ &$=$& $\mathcal{S}(M,\mathbf{P}_{p-1},\mathbf{C}_{p-1},\mathbf{K}_{p-1},\mathbf{r}_{p-1},\mathbf{u}_{p-1,n+1,0},\mathbf{u}_{p-1})$ & \rdelim\}{5}{0mm}[{\rotatebox[origin=c]{-90}{Prolong.}}]\\
					$\mathbf{\Delta}_{p-1}$ &$=$& $\mathbf{u}_{p-1,n+1,0} - \mathbf{u}_{p-1,n+1,M}$ \\
					$\mathbf{\Delta}_{p}$ &$=$& $\pmb{\pi}_p\mathbf{\Delta}_{p-1}$ \\
					$\mathbf{v}_{p,n+1,0}$ &$=$& $\mathbf{u}_{p,n+1,M} - \mathbf{\Delta}_p$ \\
					$\mathbf{u}_{p,n+1}$ &$=$& $\mathbf{R}_{p,M}\mathbf{v}_{p,n+1,0}$ \\
					\bottomrule \\
			\end{tabular}
		\end{table}
		\endgroup
		
		\begingroup
		\setlength{\tabcolsep}{2pt}
		\begin{table}[tbhp]
		    \overfullrule=0pt
		    \centering
			\caption{\label{tab:genv_steps}General $V$-cycle steps.}
			\begin{tabular}{rcll}
				for $l\in\{p, \dots, (l_\mathrm{min}+1)\}$: & & & \\ \toprule 
				    $\mathbf{u}_{l,n+1,M}$ &$=$& $\mathcal{S}(M,\mathbf{P}_{l},\mathbf{C}_{l},\mathbf{K}_{l},\mathbf{r}_{l},\mathbf{u}_{l,n+1,0},\mathbf{u}_{l})$ & \rdelim\}{5}{0mm}[{\rotatebox[origin=c]{-90}{Restriction}}]\\
					$\mathbf{d}_{l}$ &$=$& $\mathbf{r}_l - (\mathbf{T}_{l}\mathbf{u}_{l,n+1,M} - C_B\mathbf{u}_{l,n+1,0})$ \\
					$\mathbf{u}_{l-1,n+1,0}$ &$=$& $\pmb{\rho}_{l-1}\mathbf{u}_{l,n+1,M}$\\
					$\mathbf{d}_{l-1}$ &$=$& $\pmb{\rho}_{l-1}\mathbf{d}_l$ \\
					$\mathbf{r}_{l-1}$ &$=$& $\mathbf{T}_{l-1}\mathbf{u}_{l-1,n+1,0} - C_B\pmb{\rho}_{l-1}\mathbf{u}_{l,n+1,0} + \mathbf{d}_{l-1}$ \\
				for $l \in\{l_\mathrm{min}, \dots, (p-1)\}$: & & & \\\midrule
				    $\mathbf{u}_{l,n+1,M}$ &$=$& $\mathcal{S}(M,\mathbf{P}_{l},\mathbf{C}_{l},\mathbf{K}_{l},\mathbf{r}_{l},\mathbf{u}_{l,n+1,0},\mathbf{u}_{l})$ & \rdelim\}{4}{0mm}[{\rotatebox[origin=c]{-90}{Prolong.}}] \\
				    $\mathbf{\Delta}_{l}$ &$=$& $\mathbf{u}_{l,n+1,0} - \mathbf{u}_{l,n+1,M}$ \\
				    $\mathbf{\Delta}_{l+1}$ &$=$& $\pmb{\pi}_{l+1}\mathbf{\Delta}_{l}$ \\
				    $\mathbf{u}_{l+1,n+1,0}$ &$=$& $\mathbf{u}_{l+1,n+1,M} - \mathbf{\Delta}_{l+1}$ \\ \bottomrule
				    $\mathbf{u}_{p,n+1}$ &$=$& $\mathcal{S}(M,\mathbf{P}_{p},\mathbf{C}_{p},\mathbf{K}_{p},\mathbf{r}_{p},\mathbf{u}_{p,n+1,0},\mathbf{u}_{p})$ & \\
			\end{tabular}
		\end{table}
		\endgroup
		
		From the procedure defined in \cref{tab:simplev_steps,tab:genv_steps}, it can be understood that all the steps may be framed as an operation on the initial solution $\mathbf{u}_{p,n+1,0}$. It is significantly simpler to treat some steps independently and pass the result rather than formulating a single operator to act on $\mathbf{u}_{p,n+1,0}$. However, to this end, all the matrix operators at each step may be written as a polynomials in terms of $\mathbf{Q}$, and in consequence, the eigenvalues of the whole system may be found if the Bloch wave is again applied
		\begin{subequations}
			\begin{align}
				\mathbf{u}_{p,n+1} =&\: \mathbf{S}\mathbf{u}_{p,n}, \\
				=&\: \exp{(-\imath k\Delta t)}\imath k\mathbf{W}\mathbf{\Lambda}_S\mathbf{W}^{-1}\mathbf{u}_{p,n}, \\
				\exp{(\imath k\Delta t)}\mathbf{W}^{-1}\mathbf{u}_{p,n+1} =&\: \imath k\mathbf{\Lambda}_S\mathbf{W}^{-1}\mathbf{u}_{p,n},
			\end{align}				
		\end{subequations}
		where $\mathbf{S}$ is the transformation of the full system. This will enable us to examine how the energy in distributed among the spatial modes as a result of $\mathbf{W}$ being constant. Furthermore, we will define the contraction factor as
		\begin{equation}
		    \gamma = \left[\frac{\|\mathbf{e}_{m+1}\|_2 - \|\mathbf{e}_{m}\|_2}{p+1}\right]^{(n_{sp}+n_{sp}^\prime)^{-1}},
		\end{equation}
		where $n_{sp}$ is the number of smoothing iterations at the finest level applied at the beginning of the cycle and $n_{sp}^\prime$ is equivalently the number of smoothing iterations at the end of the cycle.

        \begin{figure}[tbhp]
	        \centering
            \subfloat[Error.]{\label{fig:frdg4_pmg_ad_err} \adjustbox{width=0.48\linewidth,valign=t}{    \begin{tikzpicture}
		\begin{axis}[name=plot1,xlabel={$\tau$},ylabel={$\|\mathbf{e}\|_2/(p+1)$},
            axis line style={latex-latex},
		    axis y line=left,
            axis x line=right,
            xtick={0,0.1,0.2,0.3,0.4,0.5},
		    xmin=0,xmax=0.5,
		    x tick label style={
        		/pgf/number format/.cd,
            	fixed,
            	fixed zerofill,
            	precision=1,
        	/tikz/.cd
    		},
            ymode=log,
    		ymin=3e-5,ymax=1e-1,
    		legend style={at={(0.98,0.95)},anchor=north east,font=\small},
    		legend cell align={right},
    		style={font=\large}
    		]
    		
			\addplot[color={Set1-A}, style={very thick}]table[x=Tn,y=Epi8,col sep=comma,unbounded coords=jump]{./Figs/data/FRDGp4_bdf2_ad.csv};
			\addlegendentry{Base}
			
			\addplot[color={Set1-B}, style={very thick}]table[x=Tn,y=Epi8,col sep=comma,unbounded coords=jump]{./Figs/data/FRDGp4_bdf2_pmg1V_ad.csv};
			\addlegendentry{\emph{V}-cycle, $n_s=1$}
			
			\addplot[color={Set1-C}, style={very thick}]table[x=Tn,y=Epi8,col sep=comma,unbounded coords=jump]{./Figs/data/FRDGp4_bdf2_pmg1W_ad.csv};
			\addlegendentry{$\mathrm{\emph{W}}_2$-cycle, $n_s=1$}
			
			\addplot[color={Set1-D}, style={very thick}]table[x=Tn,y=Epi8,col sep=comma,unbounded coords=jump]{./Figs/data/FRDGp4_bdf2_pmg3V_ad.csv};
			\addlegendentry{\emph{V}-cycle, $n_s=3$}
			
			\addplot[color={Set1-E}, style={very thick}]table[x=Tn,y=Epi8,col sep=comma,unbounded coords=jump]{./Figs/data/FRDGp4_bdf2_pmg1Vap3_ad.csv};
			\addlegendentry{$\mathrm{\emph{V}}_\mathrm{AP}$-cycle, $n_s=1$}
			
			
			
			
			\addplot[color={Set1-A}, style={very thick,densely dashed}]table[x=Tn,y=Epi16,col sep=comma,unbounded coords=jump]{./Figs/data/FRDGp4_bdf2_ad.csv};
			
			\addplot[color={Set1-B}, style={very thick,densely dashed}]table[x=Tn,y=Epi16,col sep=comma,unbounded coords=jump]{./Figs/data/FRDGp4_bdf2_pmg1V_ad.csv};
			
			\addplot[color={Set1-C}, style={very thick,densely dashed}]table[x=Tn,y=Epi16,col sep=comma,unbounded coords=jump]{./Figs/data/FRDGp4_bdf2_pmg1W_ad.csv};
			
			\addplot[color={Set1-D}, style={very thick,densely dashed}]table[x=Tn,y=Epi16,col sep=comma,unbounded coords=jump]{./Figs/data/FRDGp4_bdf2_pmg3V_ad.csv};
			
			\addplot[color={Set1-E}, style={very thick,densely dashed}]table[x=Tn,y=Epi16,col sep=comma,unbounded coords=jump]{./Figs/data/FRDGp4_bdf2_pmg1Vap3_ad.csv};
			
		\end{axis} 		
	\end{tikzpicture}}}
            ~
            \subfloat[Contraction.]{\label{fig:frdg4_pmg_ad_con} \adjustbox{width=0.48\linewidth,valign=t}{    \begin{tikzpicture}
		\begin{axis}[name=plot1,xlabel={$\tau$},ylabel={$\gamma$},
		    axis line style={latex-latex},
            axis y line=left,
            axis x line=right,
            xtick={0,0.1,0.2,0.3,0.4,0.5},
		    xmin=0,xmax=0.5,
    		x tick label style={
        		/pgf/number format/.cd,
            	fixed,
            	fixed zerofill,
            	precision=1,
        	/tikz/.cd
    		},
		    ylabel style={rotate={-90}},
    		ymin=0.73,ymax=1.01,
    		y tick label style={
        		/pgf/number format/.cd,
            	fixed,
            	fixed zerofill,
            	precision=2,
        	/tikz/.cd
    		},
    		legend style={at={(0.98,0.1)},anchor=south east,font=\small},
    		legend cell align={right},
    		style={font=\large}
    		]
    		
			\addplot[color={Set1-A}, style={very thick,densely dashed}]table[x=Tn,y=dEpi16,col sep=comma,unbounded coords=jump]{./Figs/data/FRDGp4_bdf2_ad.csv};
			\addlegendentry{Base}
			
			\addplot[color={Set1-B}, style={very thick,densely dashed}]table[x=Tn,y=dEpi16,col sep=comma,unbounded coords=jump]{./Figs/data/FRDGp4_bdf2_pmg1V_ad.csv};
			\addlegendentry{\emph{V}-cycle, $n_s=1$}
			
			\addplot[color={Set1-C}, style={very thick,densely dashed}]table[x=Tn,y=dEpi16,col sep=comma,unbounded coords=jump]{./Figs/data/FRDGp4_bdf2_pmg1W_ad.csv};
			\addlegendentry{$\mathrm{\emph{W}}_2$-cycle, $n_s=1$}
			
			\addplot[color={Set1-D}, style={very thick,densely dashed}]table[x=Tn,y=dEpi16,col sep=comma,unbounded coords=jump]{./Figs/data/FRDGp4_bdf2_pmg3V_ad.csv};
			\addlegendentry{\emph{V}-cycle, $n_s=3$}
			
			\addplot[color={Set1-E}, style={very thick,densely dashed}]table[x=Tn,y=dEpi16,col sep=comma,unbounded coords=jump]{./Figs/data/FRDGp4_bdf2_pmg1Vap3_ad.csv};
			\addlegendentry{$\mathrm{\emph{V}}_\mathrm{AP}$-cycle, $n_s=1$}
			
			
			

		\end{axis} 		
	\end{tikzpicture}}}
            \caption{\label{fig:frdg4_pmgV_ad_comp}Error comparison of dual time FRDG, $p=4$, with SSPRK3 and BDF2 for constant $\Delta t/\Delta\tau=10$, $\Delta\tau=\num{7e-3}$, $\mu=0.5$, and $\alpha=(1,0.5)$. For $\hat{k}=\pi/8$ (\emph{solid}), and $\pi/16$ (\emph{dashed}).}
        \end{figure}
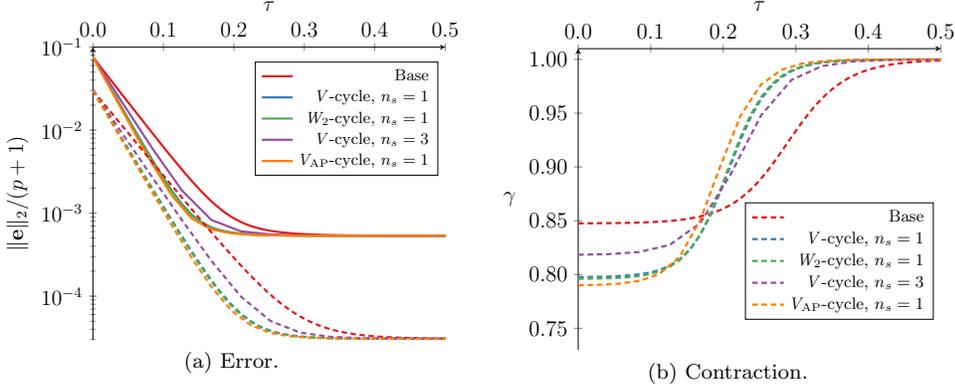
        
        \cref{fig:frdg4_pmgV_ad_comp} exemplifies the effect of \emph{p}-multigrid on the convergence of the dual time scheme. Here, advection-diffusion was considered for several cycles where $n_s$ is the number of SSPRK3 smoothing steps per level. The \emph{v}-cycle with subscript AP has additional prolongation smoothing steps with the prolongation smoothing set to three. An example cycle is given in \cref{fig:ap_v}. The pseudo time shown for the \emph{p}-multigrid cases is the cumulative time at the finest \emph{p} level, i.e., $\tau=(n_{sp}+n_{sp}^\prime)n_\mathrm{cycle}\Delta\tau$. In all cases, \emph{p}-multigrid increased the rate of convergence, however, it is clear that fewer smoothing steps during the restriction portion of the cycle was beneficial to convergence. A corollary observation is that making the \emph{V}-cycle asymmetric with addition prolongation smoothing could further increase convergence. This is due to the larger differences in pseudo time between levels causing larger deficit source terms which, upon prolongation, lead to the need for more smoothing steps so that they are adequately relaxed into the solution. From the \emph{V}-cycle with $n_s=3$ in \cref{fig:frdg4_pmgV_ad_comp}, it is clear that the prolongation smoothing requirement does not grow linearly with the restriction smoothing otherwise we would expect to see results closer to the $n_s=1$ case. 
        
        It was also observed in all cases that the number of overall iterations to converge is limited by the lower wavenumbers. This result may be expected and can be understood from the longer half-lives of the these waves due to the lower dissipation when considered in the fully discrete form. As the viscosity was decreased, the effectiveness of \emph{p}-multigrid  was found to decrease. However, additional prolongation was still found to be effective.
        
        \begin{figure}[tbhp]
	        \centering
            \subfloat[Primary mode.]{\label{fig:frdg4_pmg_ad_beta} \adjustbox{width=0.48\linewidth,valign=b}{    \begin{tikzpicture}
		\begin{axis}[name=plot1,xlabel={$n_c$},ylabel={$|\pmb{\beta}_p|$},
		    axis line style={latex-latex},
            axis y line=left,
            axis x line=left,
		    xmin=0,xmax=10,
    		ymin=2.15,ymax=2.25,
    		y tick label style={
        		/pgf/number format/.cd,
            	fixed,
            	fixed zerofill,
            	precision=2,
        	/tikz/.cd
    		},
    		legend style={at={(0.95,0.95)},anchor=north east,font=\small},
    		legend cell align={right},
    		style={font=\large}
    		]
    		
			\addplot[color={Set1-A}, style={very thick}]table[x expr=\thisrowno{0}/0.007,y=Beta_4,col sep=comma,unbounded coords=jump]{./Figs/data/FRDGp4_BDF2_ad_re2_beta.csv};
			\addlegendentry{Base}
			
			\addplot[color={Set1-B}, style={very thick}]table[x expr=\thisrowno{0}/0.007,y=Beta_1V_4,col sep=comma,unbounded coords=jump]{./Figs/data/FRDGp4_BDF2_ad_re2_beta.csv};
			\addlegendentry{\emph{V}-cycle, $n_s=1$}
			
			\addplot[color={Set1-C}, style={very thick}]table[x expr=\thisrowno{0}/0.007,y=Beta_3V_4,col sep=comma,unbounded coords=jump]{./Figs/data/FRDGp4_BDF2_ad_re2_beta.csv};
		    \addlegendentry{\emph{V}-cycle, $n_s=3$}
			
			\addplot[color={Set1-D}, style={very thick}]table[x expr=\thisrowno{0}/0.007,y=Beta_1Vap3_4,col sep=comma,unbounded coords=jump]{./Figs/data/FRDGp4_BDF2_ad_re2_beta.csv};
			\addlegendentry{$\mathrm{\emph{V}}_\mathrm{AP}$-cycle, $n_s=1$}
			
		\end{axis} 		
	\end{tikzpicture}}}
            ~
            \subfloat[Secondary mode.]{\label{fig:frdg4_pmg_ad_beta2} \adjustbox{width=0.48\linewidth,valign=b}{    \begin{tikzpicture}
		\begin{axis}[name=plot1,xlabel={$n_c$},ylabel={$|\pmb{\beta}_p|$},
		    axis line style={latex-latex},
            axis y line=left,
            axis x line=left,
		    xmin=0,xmax=10,
    		ymin=0,
    		ymax=6.5e-4,
    		ytick={0,1e-4,2e-4,3e-4,4e-4,5e-4,6e-4},
    		legend style={at={(0.98,0.98)},anchor=north east,font=\small},
    		legend cell align={right},
    		style={font=\large}
    		]
    		
			\addplot[color={Set1-A}, style={very thick}]table[x expr=\thisrowno{0}/0.007,y=Beta_5,col sep=comma,unbounded coords=jump]{./Figs/data/FRDGp4_BDF2_ad_re2_beta.csv};
			\addlegendentry{Base}
			
			\addplot[color={Set1-B}, style={very thick}]table[x expr=\thisrowno{0}/0.007,y=Beta_1V_5,col sep=comma,unbounded coords=jump]{./Figs/data/FRDGp4_BDF2_ad_re2_beta.csv};
			\addlegendentry{\emph{V}-cycle, $n_s=1$}
			
			\addplot[color={Set1-C}, style={very thick}]table[x expr=\thisrowno{0}/0.007,y=Beta_3V_5,col sep=comma,unbounded coords=jump]{./Figs/data/FRDGp4_BDF2_ad_re2_beta.csv};
			\addlegendentry{\emph{V}-cycle, $n_s=3$}
			
			\addplot[color={Set1-D}, style={very thick}]table[x expr=\thisrowno{0}/0.007,y=Beta_1Vap3_5,col sep=comma,unbounded coords=jump]{./Figs/data/FRDGp4_BDF2_ad_re2_beta.csv};
			\addlegendentry{$\mathrm{\emph{V}}_\mathrm{AP}$-cycle, $n_s=1$}
			
		\end{axis} 		
	\end{tikzpicture}}}
            \caption{\label{fig:frdg4_beta_ad_comp}Primary and secondary mode energy at $\hat{k}=\pi/16$ for FRDG $p=4$, BDF2, and SSPRK3, with $\Delta\tau=0.007$, $\Delta t/\Delta\tau=10$, $\mu=0.5$, and $\alpha=(1,0.5)$. Note the change in \emph{y}-axis scaling between figures.}
        \end{figure}
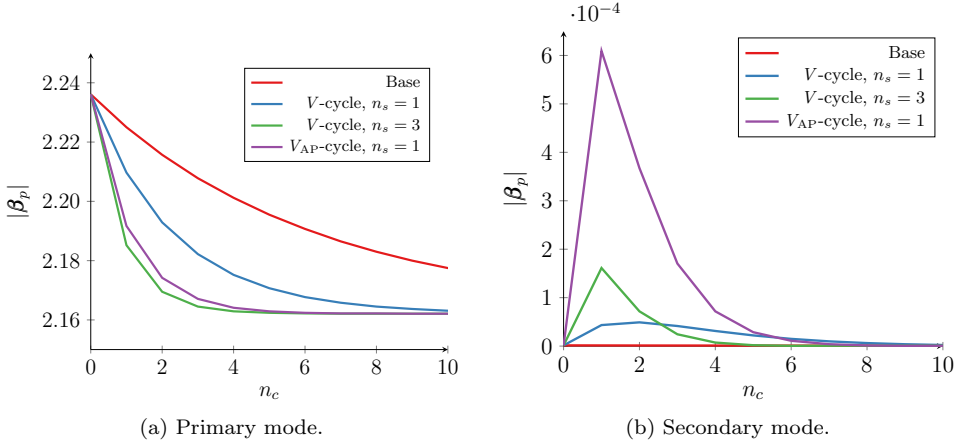
        
        Further insight as to why the additional prolongation \emph{V}-cycles have improved contraction rates can be gained from inspection of the way in which the solution energy is distributed among the modes of the spatial system. \cref{fig:frdg4_beta_ad_comp} shows the energy in the primary and secondary modes for several cycle configurations. The additional prolongation steps, in both cases, causes a greater redistribution of energy from the primary mode to the secondary mode. When considered with the knowledge that the secondary modes have shorter half-lives~\cite{Trojak2018}, the mechanism of convergence acceleration is understood to come from this redistribution. The additional restriction smoothing steps in the $n_s=3$ case diminishes the redistribution, and hence is why the contraction factor in \cref{fig:frdg4_pmg_ad_con} show poorer acceleration. If additional restriction smoothing alone is considered, then the effect on redistribution compared to the $n_s=1$ case is negligible, which is concurrent with redistribution being due to the prolongation correction as may have been anticipated.
        
        \begin{figure}[tbhp]
	        \centering
            \subfloat[BDF2.]
            {\label{fig:FRDG_bdf2_cont}\adjustbox{width=0.48\linewidth,valign=b}{\begin{tikzpicture}
		\begin{axis}[name=plot1,xlabel={$\Delta t/\Delta\tau$},ylabel={$\gamma$},
            axis line style={latex-latex},
		    axis x line=left,
            axis y line=left,
            xmode=log,
		    xmin=1,
		    xmax=1e2,
		    ylabel style={rotate={-90}},
    		ymin=0,
    		ymax=1,
    		y tick label style={
        		/pgf/number format/.cd,
            	fixed,
            	fixed zerofill,
            	precision=1,
        	/tikz/.cd
    		},
    		legend style={at={(0.95,0.05)},anchor=south east,font=\small},
    		legend cell align={left},
    		style={font=\large}
    		]

            \addplot[color={Set1-A}, style={very thick}]
                table[x=Beta,y=Gb,col sep=comma,unbounded coords=jump]{./Figs/data/FRDGp4_BDF2_ad_re10_C.csv};
			\addlegendentry{Base}
			
            \addplot[color={Set1-B}, style={very thick}]
                table[x=Beta,y=G1Vap3,col sep=comma,unbounded coords=jump]{./Figs/data/FRDGp4_BDF2_ad_re10_C.csv};
			\addlegendentry{$\mathrm{\emph{V}}_\mathrm{AP}$-cycle}
			
			 \addplot[color=black,style={very thick}] coordinates{(5.0941,0) (5.0941,1)};
			\addlegendentry{Max. decrease}

            \addplot[color={Set1-A}, style={very thick,densely dashed}]
                table[x=Beta,y=Gb,col sep=comma,unbounded coords=jump]{./Figs/data/FRDGp3_BDF2_ad_re10_C.csv};
			
            \addplot[color={Set1-B}, style={very thick,densely dashed}]
                table[x=Beta,y=G1Vap3,col sep=comma,unbounded coords=jump]{./Figs/data/FRDGp3_BDF2_ad_re10_C.csv};
                
			 \addplot[color=black,style={very thick,densely dashed}] coordinates{(3.5112,0) (3.5112,1)};

    \end{axis}    		
\end{tikzpicture}}}
            ~
            \subfloat[BDF3.]
            {\label{fig:FRDG_bdf3_cont}\adjustbox{width=0.48\linewidth,valign=b}{\begin{tikzpicture}
		\begin{axis}[name=plot1,xlabel={$\Delta t/\Delta\tau$},ylabel={$\gamma$},
            axis line style={latex-latex},
		    axis x line=left,
            axis y line=left,
            xmode=log,
		    xmin=1,
		    xmax=1e2,
		    ylabel style={rotate={-90}},
    		ymin=0,
    		ymax=1,
    		y tick label style={
        		/pgf/number format/.cd,
            	fixed,
            	fixed zerofill,
            	precision=1,
        	/tikz/.cd
    		},
    		legend style={at={(0.95,0.05)},anchor=south east,font=\small},
    		legend cell align={left},
    		style={font=\large}
    		]

            \addplot[color={Set1-A}, style={very thick}]
                table[x=Beta,y=Gb,col sep=comma,unbounded coords=jump]{./Figs/data/FRDGp4_BDF3_ad_re10_C.csv};
			\addlegendentry{Base}
			
            \addplot[color={Set1-B}, style={very thick}]
                table[x=Beta,y=G1Vap3,col sep=comma,unbounded coords=jump]{./Figs/data/FRDGp4_BDF3_ad_re10_C.csv};
			\addlegendentry{$\mathrm{\emph{V}}_\mathrm{AP}$-cycle}
			
			 \addplot[color=black,style={very thick}] coordinates{(6.1359,0) (6.1359,1)};
			\addlegendentry{Max. decrease}
			
            \addplot[color={Set1-A}, style={very thick,densely dashed}]
                table[x=Beta,y=Gb,col sep=comma,unbounded coords=jump]{./Figs/data/FRDGp3_BDF3_ad_re10_C.csv};
			
            \addplot[color={Set1-B}, style={very thick,densely dashed}]
                table[x=Beta,y=G1Vap3,col sep=comma,unbounded coords=jump]{./Figs/data/FRDGp3_BDF3_ad_re10_C.csv};
			
			 \addplot[color=black,style={very thick,densely dashed}] coordinates{(4.2292,0) (4.2292,1)};

    \end{axis}    		
\end{tikzpicture}}}
            \caption{\label{fig:FRDG_cont}Initial contraction factor for FRDG with BDF and SSPRK3 at $k=(p+1)\pi/16$, $\Delta\tau=0.078\Delta t_\mathrm{max,A}$, $\mu=0.1$, and $\alpha=(1,0.5)$. With spatial orders $p=4$ (\emph{solid}) and $p=3$ (\emph{dashed}).}
        \end{figure}
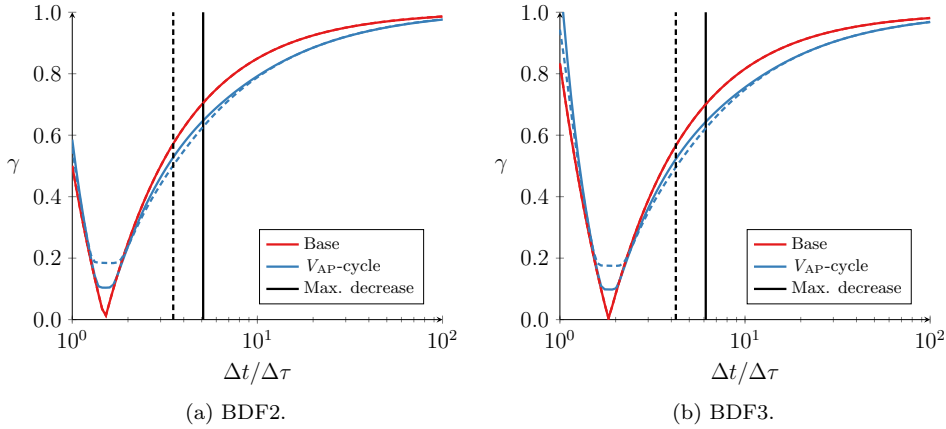
        
        Subsequently, for a constant wavenumber, values of $\Delta t/\Delta\tau$ were swept through and the contraction was found, the results of which are presented in \cref{fig:FRDG_cont}. Here, we have only used a \emph{V}-cycle with additional prolongation as this offered the best performance. From this data the diminishing returns of using multigrid to accelerate dual-time for large ratios of $\Delta t/\Delta\tau$ is seen. This is due to the large time scales in the hyperbolic component  of the system becoming dominant, therefore the dual-time convergence is primarily just dependent on the number of iterations. We have marked the points on each diagram where the ratio of contraction between the base scheme and \emph{p}-multigrid is largest. A move from BDF2 to BDF3, for both spatial order tested, resulted in the maximal point increasing by ${\sim}20\%$. 
        
        As a point of comparison, the element Jacobi method coupled to BDF was also considered, both with and without \emph{p}-multigrid acceleration. A brief description, and associated definitions of this technique are included in \cref{app:ej}, and after the definition of the EJ matrix, the earlier derivations may be followed to apply the \emph{p}-multigrid methodology. 
        
        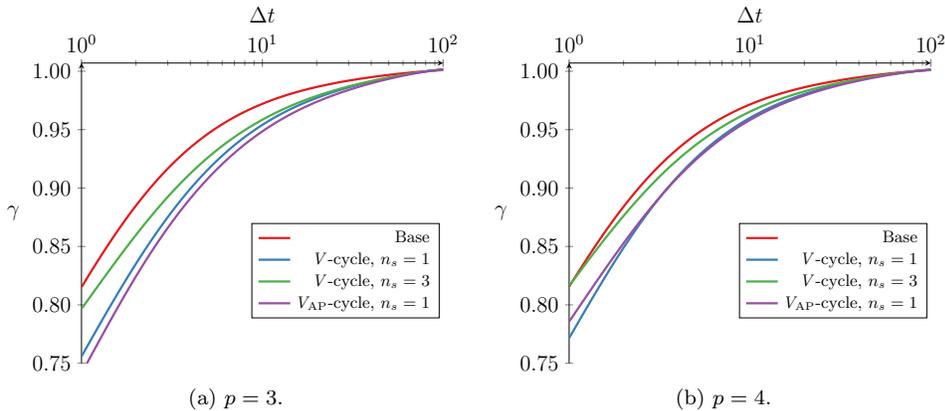
\begin{figure}[tbhp]
	        \centering
            \subfloat[$p=3$.]
            {\label{fig:FRDG3_ej_cont}\adjustbox{width=0.48\linewidth,valign=t}{\begin{tikzpicture}
		\begin{axis}[name=plot1,xlabel={$\Delta t$},ylabel={$\gamma$},
            axis line style={latex-latex},
            axis y line=left,
            axis x line=right,
            xmode=log,
		    xmin=1,
		    xmax=1e2,
		    ylabel style={rotate={-90}},
    		ymin=0.75,
    		ymax=1.007,
    		y tick label style={
        		/pgf/number format/.cd,
            	fixed,
            	fixed zerofill,
            	precision=2,
        	/tikz/.cd
    		},
    		legend style={at={(0.99,0.15)},anchor=south east,font=\small},
    		legend cell align={right},
    		style={font=\large}
    		]

            \addplot[color={Set1-A}, style={very thick}]
                table[x=dT,y=Gb,col sep=comma,unbounded coords=jump]{./Figs/data/FRDGp3_BDF2_ej_ad_re10_C.csv};
			\addlegendentry{Base}
			
            \addplot[color={Set1-B}, style={very thick}]
                table[x=dT,y=G1V,col sep=comma,unbounded coords=jump]{./Figs/data/FRDGp3_BDF2_ej_ad_re10_C.csv};
			\addlegendentry{\emph{V}-cycle, $n_s=1$}
			
            \addplot[color={Set1-C}, style={very thick}]
                table[x=dT,y=G3V,col sep=comma,unbounded coords=jump]{./Figs/data/FRDGp3_BDF2_ej_ad_re10_C.csv};
			\addlegendentry{\emph{V}-cycle, $n_s=3$}
			
            \addplot[color={Set1-D}, style={very thick}]
                table[x=dT,y=G1Vap3,col sep=comma,unbounded coords=jump]{./Figs/data/FRDGp3_BDF2_ej_ad_re10_C.csv};
			\addlegendentry{$\mathrm{\emph{V}}_\mathrm{AP}$-cycle, $n_s=1$}

    \end{axis}    		
\end{tikzpicture}}}
            ~
            \subfloat[$p=4$.]
            {\label{fig:FRDG4_ej_cont}\adjustbox{width=0.48\linewidth,valign=t}{\begin{tikzpicture}
		\begin{axis}[name=plot1,xlabel={$\Delta t$},ylabel={$\gamma$},
            axis line style={latex-latex},
            axis y line=left,
            axis x line=right,
            xmode=log,
		    xmin=1,
		    xmax=1e2,
		    ylabel style={rotate={-90}},
    		ymin=0.75,
    		ymax=1.007,
    		y tick label style={
        		/pgf/number format/.cd,
            	fixed,
            	fixed zerofill,
            	precision=2,
        	/tikz/.cd
    		},
    		legend style={at={(0.99,0.15)},anchor=south east,font=\small},
    		legend cell align={right},
    		style={font=\large}
    		]

            \addplot[color={Set1-A}, style={very thick}]
                table[x=dT,y=Gb,col sep=comma,unbounded coords=jump]{./Figs/data/FRDGp4_BDF2_ej_ad_re10_C.csv};
			\addlegendentry{Base}
			
            \addplot[color={Set1-B}, style={very thick}]
                table[x=dT,y=G1V,col sep=comma,unbounded coords=jump]{./Figs/data/FRDGp4_BDF2_ej_ad_re10_C.csv};
			\addlegendentry{\emph{V}-cycle, $n_s=1$}
			
            \addplot[color={Set1-C}, style={very thick}]
                table[x=dT,y=G3V,col sep=comma,unbounded coords=jump]{./Figs/data/FRDGp4_BDF2_ej_ad_re10_C.csv};
			\addlegendentry{\emph{V}-cycle, $n_s=3$}
			
            \addplot[color={Set1-D}, style={very thick}]
                table[x=dT,y=G1Vap3,col sep=comma,unbounded coords=jump]{./Figs/data/FRDGp4_BDF2_ej_ad_re10_C.csv};
			\addlegendentry{$\mathrm{\emph{V}}_\mathrm{AP}$-cycle, $n_s=1$}

    \end{axis}    		
\end{tikzpicture}}}
            \caption{\label{fig:FRDG_ej_cont}Initial contraction factor for FRDG with BDF2 and EJ at $k=\pi/100$, $\kappa=0.5$, $\mu=0.1$, and $\alpha=(1,0.5)$.}
        \end{figure}
        
        The contraction factor for the element Jacobi method is presented in \cref{fig:FRDG_ej_cont} for BDF2 at two spatial orders. Similar trends to those observed for the dual-time scheme are seen here, with additional prolongation smoothing being favourable. However, at higher spatial orders and lower time steps, additional prolongation and the $n_s=3$ cycle saw a reduction in their benefit. As the $n_s=1$ cycle maintained the improved contraction, this degradation is due to one smoothing step being sufficient in this less stiff range of $\Delta t$.

    \subsection{\emph{p}-Multigrid Acceleration}
        As has been confirmed here, \emph{p}-multigrid does not have as greater benefit to accelerate the convergence of the coupled ERK-BDF dual-time system for hyperbolic equations. This is evident when considering the contraction factor in \cref{fig:FRDG_cont} in the limit as $\Delta t/\Delta\tau\rightarrow\infty$ where the hyperbolic time scales become dominant, and with \emph{p}-multigrid providing a greater degree of acceleration for elliptic-hyperbolic equations. As has been discussed in the literature, this is due to the local dependency of hyperbolic equations compared to the global dependency of elliptic problems~\cite{Katzer1991}, and it follows that the convergence of hyperbolic components here are dependent on the convection time of waves in the system. This is not to say that \emph{p}-multigrid cannot be effective for hyperbolic problems. For example, it can be effective when employing Newton--Krylov approaches with large time steps as in the limit the system becomes elliptic.
        
        One method to further accelerate the dual-time \emph{p}-multigrid investigated here is a procedure where the pseudo-time step was increased at coarser \emph{p}-multigrid levels, see Loppi~\etal~\cite{Loppi2019}. In this method, a factor was introduced such that the pseudo-time step is defined as
        \begin{equation}
            \Delta\tau_i = \Delta\tau(f_\tau^{p-i}),
        \end{equation}
        where $\Delta\tau_i$ is the pseudo-time step at the degree $i$ \emph{p}-multigrid level. When setting $f_\tau$, care must be taken such that the CFL limits imposed through \cref{fig:FRDG_ad_cfl,fig:FRDG_BDF_CFL} are not exceeded. 
        
        This method will allow more rapid advection---as well as diffusion---of waves at the coarser levels. Implicitly these waves are of lower frequency and consequently are more challenging to converge due to their large length and time scales. This may also pose a problem if the corrections are not sufficiently relaxed into the finer multigrid levels as the corrections are likely to be large due to the different pseudo time steps, allowing error to accumulate in the solution.
        
        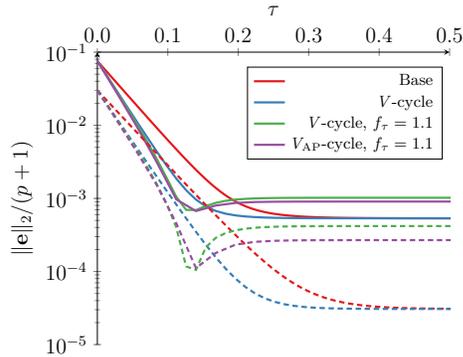
\begin{figure}[tbhp]
	        \centering
            \adjustbox{width=0.48\linewidth,valign=b}{    \begin{tikzpicture}
		\begin{axis}[name=plot1,xlabel={$\tau$},ylabel={$\|\mathbf{e}\|_2/(p+1)$},
            axis line style={latex-latex},
            axis y line=left,
            axis x line=right,
            xtick={0,0.1,0.2,0.3,0.4,0.5},
		    xmin=0,xmax=0.5,
		    x tick label style={
        		/pgf/number format/.cd,
            	fixed,
            	fixed zerofill,
            	precision=1,
        	/tikz/.cd
    		},
            ymode=log,
    		ymin=1e-5,ymax=1e-1,
    		legend style={at={(0.98,0.95)},anchor=north east,font=\small},
    		legend cell align={right},
    		style={font=\large}
    		]
    		
			\addplot[color={Set1-A}, style={very thick}]table[x=Tn,y=Epi8,col sep=comma,unbounded coords=jump]{./Figs/data/FRDGp4_bdf2_ad.csv};
			\addlegendentry{Base}
			
			\addplot[color={Set1-B}, style={very thick}]table[x=Tn,y=Epi8,col sep=comma,unbounded coords=jump]{./Figs/data/FRDGp4_bdf2_pmg1V_ad.csv};
			\addlegendentry{\emph{V}-cycle}
			
			\addplot[color={Set1-C}, style={very thick}]table[x=Tn,y=Epi8,col sep=comma,unbounded coords=jump]{./Figs/data/FRDGp4_bdf2_pmg1V_ad_fc11.csv};
			\addlegendentry{\emph{V}-cycle, $f_\tau=1.1$}
			
			\addplot[color={Set1-D}, style={very thick}]table[x=Tn,y=Epi8,col sep=comma,unbounded coords=jump]{./Figs/data/FRDGp4_bdf2_pmg1Vap3_ad_fc11.csv};
			\addlegendentry{$\mathrm{\emph{V}}_\mathrm{AP}$-cycle, $f_\tau=1.1$}
			
			\addplot[color={Set1-A}, style={very thick,densely dashed}]table[x=Tn,y=Epi16,col sep=comma,unbounded coords=jump]{./Figs/data/FRDGp4_bdf2_ad.csv};
			
			\addplot[color={Set1-B}, style={very thick,densely dashed}]table[x=Tn,y=Epi16,col sep=comma,unbounded coords=jump]{./Figs/data/FRDGp4_bdf2_pmg1V_ad.csv};
			
			\addplot[color={Set1-C}, style={very thick,densely dashed}]table[x=Tn,y=Epi16,col sep=comma,unbounded coords=jump]{./Figs/data/FRDGp4_bdf2_pmg1V_ad_fc11.csv};
			
			\addplot[color={Set1-D}, style={very thick,densely dashed}]table[x=Tn,y=Epi16,col sep=comma,unbounded coords=jump]{./Figs/data/FRDGp4_bdf2_pmg1Vap3_ad_fc11.csv};

		\end{axis} 		
	\end{tikzpicture}}
            \caption{\label{fig:fc_pmg}Error comparison of dual time FRDG with $f_\tau$, $p=4$, with SSPRK3 and BDF2 for constant $\Delta t/\Delta\tau=10$, $\Delta\tau=\num{7e-3}$, $\mu=0.5$, and $\alpha=(1,0.5)$. For $\hat{k}=\pi/8$ (\emph{solid}), and $\pi/16$ (\emph{dashed}).}
        \end{figure}
       
        \cref{fig:fc_pmg} presents the results of applying $f_\tau=1.1$ to the \emph{p}-multigrid cycle with one smoothing step per stage. From this data, it may be concluded that rate of convergence is further increased by $f_\tau$. However, as was hypothesised, insufficient relaxation during prolongation causes the build up of error. This may be mollified by additional prolongation, with the result here using three smoothing steps per prolongation stage, but significant steady state error is still present. 
        
\section{Numerical Experiments}\label{sec:numeric}
    In order to test the analytic hypothesis about the utility of asymmetric \emph{V}-cycles, we will consider the incompressible Navier--Stokes equations solved via ACM. The governing equations in two dimensions takes the form
    \begin{subequations}
        \begin{align}
            \px{P}{\tau} + \px{(\zeta u)}{x} + \px{(\zeta v)}{y} &= 0, \\
            \px{u}{\tau} + \px{u}{t} + \px{(u^2 + P)}{x} + \px{uv}{y} &= \nu\left(\pxi{2}{u}{x} + \pxi{2}{u}{y}\right), \\
            \px{v}{\tau} + \px{v}{t} + \px{(v^2 +P)}{y} + \px{uv}{x} &= \nu\left(\pxi{2}{v}{x} + \pxi{2}{v}{y}\right),
        \end{align}
    \end{subequations}
    where $P$ is the pressure, $u$ and $v$ are the \emph{x} and \emph{y} components of velocity, respectively, $\nu$ is the kinematic viscosity, and $\zeta$ is the artificial compressibility coefficient. The numerical experiments were performed using the high-order FR solver PyFR~\cite{Witherden2014,Loppi2018}. DG recovering correction functions were used together with BR1~\cite{Bassi1997} viscous and HLLC~\cite{Elsworth1992} inviscid ACM approximate common interface flux calculations. The solution and flux points were positioned using Gauss--Legendre and Williams--Shunn~\cite{Williams2014} points for quadrilaterals and triangles, respectively.
    
    \begin{figure}[tbhp]
        \centering
        \subfloat[View of unstructured mesh.]{\label{fig:naca_mesh}\includegraphics[width=0.49\linewidth]{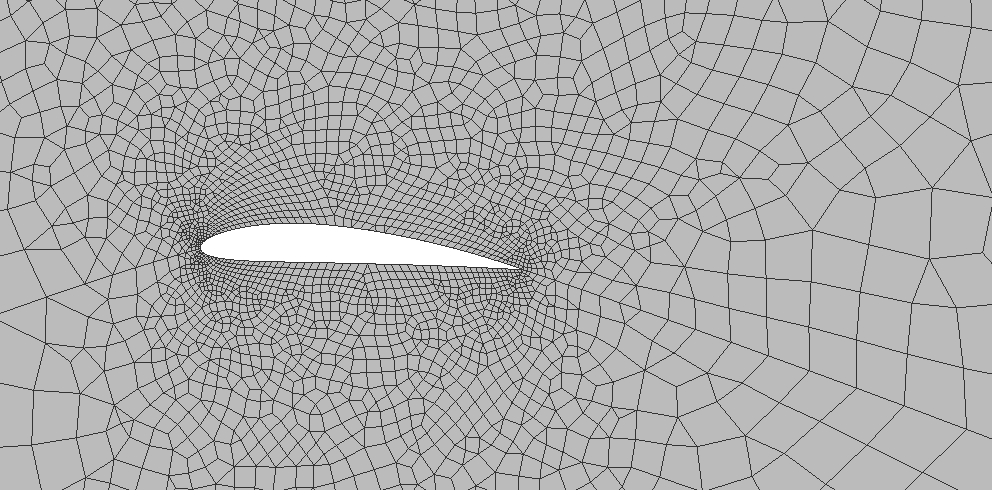}}
        ~
        \subfloat[Pressure at $t=\SI{70}{\second}$.]{\label{fig:naca_pressure}\includegraphics[width=0.49\linewidth]{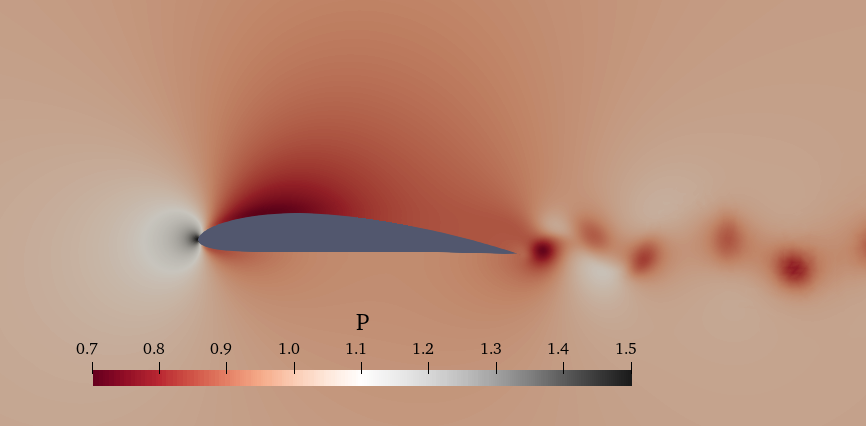}}
        \caption{NACA-4412 at $\mathrm{AoA}=2.5^\circ$ and $Re=\num{5e3}$.}
    \end{figure}
    
   The operating condition examined throughout the experiments was at an angle-of-attack~(AoA) of $2.5^\circ$, $Re=\num{5e3}$, a spatial order of $p=3$, and $\zeta=2.5$. The far-field pressure and velocity magnitude were $P_\infty=1$ and $V_\infty=1$, respectively.  A view of the mesh used can be seen in \cref{fig:naca_mesh} and is comprised of ${\sim}900$ triangles and ${\sim}3800$ quadrilaterals. Although this is a simple geometry at low Reynolds number, we chose to use a fully unstructured mesh as this better represents the typical use case for this method.
   
   For the temporal integration, BDF2 with SSPRK3 for the pseudo-time stepping was used, with $\Delta t=\num{5e-4}$, $\Delta t/\Delta\tau=5$, and $f_\tau=1.75$. A fixed number of pseudo steps per iteration of ten was used. A higher number would typically be needed for engineering a calculation; however, this was deemed to be sufficient to demonstrate the convergence acceleration in this case. The simulations were run for $75$ flows over chord and the pressure distribution at $t=\SI{75}{\second}$ is shown in \cref{fig:naca_pressure} where the vortex shedding is clearly visible. 
    
    \begin{figure}[tbhp]
        \centering
        \subfloat[Pressure.]{\label{fig:naca_conv_p} \resizebox{0.47\linewidth}{!}{    \begin{tikzpicture}
		\begin{axis}[name=plot1,xlabel={$n_c$},ylabel={$\overline{R_{n_c}}/\overline{R_1}$},
		    axis line style={latex-latex},
            axis y line=left,
            axis x line=left,
		    xmin=1,xmax=10,
    		ymin=0.80,ymax=1.007,
    		y tick label style={
        		/pgf/number format/.cd,
            	fixed,
            	fixed zerofill,
            	precision=2,
        	/tikz/.cd
    		},
    		legend style={at={(0.98,0.945)},anchor=north east,font=\small},
    		legend cell align={right},
    		style={font=\large}
    		]
    		
			\addplot[color={Set1-A}, style={very thick}]table[x=nc,y=base,col sep=comma,unbounded coords=jump]{./Figs/data/naca4412_resid_conv_p.csv};
			\addlegendentry{Base}
			
			\addplot[color={Set1-B}, style={very thick}]table[x=nc,y=1V,col sep=comma,unbounded coords=jump]{./Figs/data/naca4412_resid_conv_p.csv};
			\addlegendentry{\emph{V}-cycle, $n_s=1$}
			
			\addplot[color={Set1-C}, style={very thick}]table[x=nc,y=3V,col sep=comma,unbounded coords=jump]{./Figs/data/naca4412_resid_conv_p.csv};
			\addlegendentry{\emph{V}-cycle, $n_s=3$}
			
			\addplot[color={Set1-D}, style={very thick}]table[x=nc,y=1Vap3,col sep=comma,unbounded coords=jump]{./Figs/data/naca4412_resid_conv_p.csv};
			\addlegendentry{$\mathrm{\emph{V}}_\mathrm{AP}$-cycle, $n_s=1$}
			
		\end{axis} 		
	\end{tikzpicture}}}
        ~
        \subfloat[\emph{x}-velocity.]{\label{fig:naca_conv_u} \resizebox{0.47\linewidth}{!}{\begin{tikzpicture}
		\begin{axis}[name=plot1,xlabel={$n_c$},ylabel={$\overline{R_{n_c}}/\overline{R_1}$},
		    axis line style={latex-latex},
            axis y line=left,
            axis x line=left,
		    xmin=1,xmax=10,
    		ymin=0,ymax=1.007,
    		y tick label style={
        		/pgf/number format/.cd,
            	fixed,
            	fixed zerofill,
            	precision=2,
        	/tikz/.cd
    		},
    		legend style={at={(0.95,0.95)},anchor=north east,font=\small},
    		legend cell align={right},
    		style={font=\large}
    		]
    		
			\addplot[color={Set1-A}, style={very thick}]table[x=nc,y=base,col sep=comma,unbounded coords=jump]{./Figs/data/naca4412_resid_conv_u.csv};
			\addlegendentry{Base}
			
			\addplot[color={Set1-B}, style={very thick}]table[x=nc,y=1V,col sep=comma,unbounded coords=jump]{./Figs/data/naca4412_resid_conv_u.csv};
			\addlegendentry{\emph{V}-cycle, $n_s=1$}
			
			\addplot[color={Set1-C}, style={very thick}]table[x=nc,y=3V,col sep=comma,unbounded coords=jump]{./Figs/data/naca4412_resid_conv_u.csv};
			\addlegendentry{\emph{V}-cycle, $n_s=3$}
			
			\addplot[color={Set1-D}, style={very thick}]table[x=nc,y=1Vap3,col sep=comma,unbounded coords=jump]{./Figs/data/naca4412_resid_conv_u.csv};
			\addlegendentry{$\mathrm{\emph{V}}_\mathrm{AP}$-cycle, $n_s=1$}
			
		\end{axis} 		
	\end{tikzpicture}}}
        \caption{\label{fig:naca_conv}Mean relative residual convergence for NACA-4412 at $p=3$ over \num{1e4} implicit steps.}
    \end{figure}
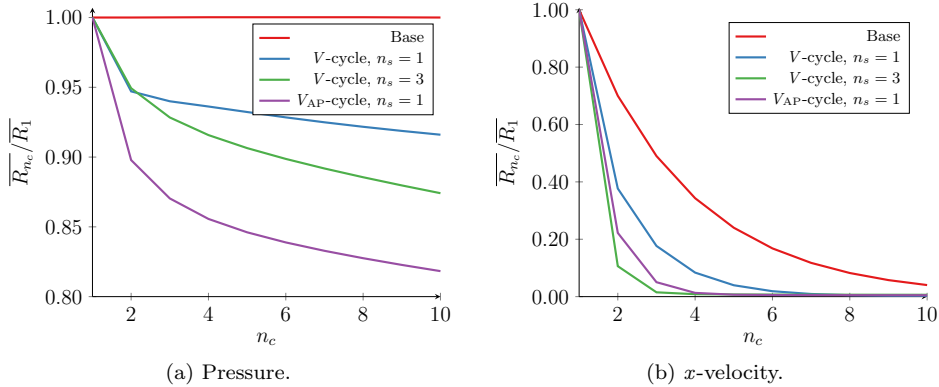
    
    To demonstrate the effect of various \emph{V}-cycles, we investigated the averaged relative residual for each cycle in dual-time. The mean residual for each cycle is normalised by the mean of the initial residual in each real time step. The results averaged for the last \num{1e4} physical-time steps, equivalent to  approximately 10 shedding cycles, is presented in \cref{fig:naca_conv}. An interesting difference in behaviour is exhibited between the pressure and velocity convergence, with pressure showing the same predicted improvement for additional prolongation, whereas for the convergence of velocity cycles, more smoothing steps caused the fastest decay in the residual. This is due to the different character of the equations; the first equation---which drives pressure---is elliptic, whereas the velocity equations are hyperbolic. Hence, the convergence of the velocity equations is chiefly a matter of advection and benefits primarily from a greater number of pseudo time iterations.
    
    The low number of pseudo-steps used here is visible for the base case from the high average pressure residual shown in \cref{tab:naca_mean_quant} and that fact that the residual factor in \cref{fig:naca_conv_p} for the base case does not show reduction. However, reduction is still seen in the velocity residual, for which the governing equation is dominantly hyperbolic and hence benefits purely from additional iterations to further convergence.
    
    \begin{table}[H]
        \centering
        \caption{\label{tab:naca_mean_quant}Metric comparison for various cycles.}
        \begin{tabular}{lcccc} \toprule
            Cycle & $n_s$ & $C_L/C_D$ & $\overline{R}_P$ &  $\overline{R}_u$ \\ \midrule
            None &  & 7.246 & \num{5.161e-2} & \num{1.450e-1} \\
            \emph{V} & $1$ & 7.072 & \num{3.547e-3} & \num{2.158e-3} \\
            \emph{V} & $3$ & 7.071 & \num{1.173e-3} & \num{2.698e-4} \\
            $\mathrm{\emph{V}}_\mathrm{AP}$ & $1$ & 7.067 & \num{1.612e-3} & \num{6.249e-04} \\ \bottomrule
        \end{tabular}
    \end{table}

\section{Conclusions}\label{sec:conclusions}
In this manuscript, we have presented a Fourier analysis of dual-time stepping with the high-order FR approach using \emph{p}-multigrid convergence acceleration.  This enables---for the first time---arbitrary multigrid cycles to be explored and analysed directly.  Employing this analysis, we have shown for the advection-diffusion equation that \emph{p}-multigrid can reduce the contraction factor by $9\%$.  Furthermore, it was also shown how performance can be improved through the use of \emph{asymmetric} cycles which contain additional prolongation steps, an observation which is supported through numerical experiments with the incompressible Navier--Stokes equations on a 2D NACA-4412.

\section*{Acknowledgements}\label{sec:ack}
We would like thank to T. Dzanic and L. Wang for aiding us in preparation of this manuscript.

\bibliographystyle{siamplain}
\bibliography{reference}

\begin{thebibliography}{10}

\bibitem{Arnone1993}
{\sc A.~Arnone, M.~Liou, and L.~Povinelli}, {\em {Multigrid Time-Accurate
  Integration of Navier--Stokes Equations}}, in 11th Computational Fluid
  Dynamics Conference, American Institute of Aeronautics and Astronautics, July
  1993, \url{https://doi.org/10.2514/6.1993-3361},
  \url{https://doi.org/10.2514/6.1993-3361}.

\bibitem{Bassi1997}
{\sc F.~Bassi and S.~Rebay}, {\em {A High-Order Accurate Discontinuous Finite
  Element Method for the Numerical Solution of the Compressible Navier--Stokes
  Equations}}, Journal of Computational Physics, 131 (1997), pp.~267--279,
  \url{https://doi.org/10.1006/jcph.1996.5572},
  \url{https://doi.org/10.1006/jcph.1996.5572}.

\bibitem{Butcher1964}
{\sc J.~C. Butcher}, {\em {On Runge--Kutta Processes of High Order}}, Journal
  of the Australian Mathematical Society, 4 (1964), pp.~179--194,
  \url{https://doi.org/10.1017/s1446788700023387},
  \url{https://doi.org/10.1017/s1446788700023387}.

\bibitem{Castonguay2013}
{\sc P.~Castonguay, D.~Williams, P.~Vincent, and A.~Jameson}, {\em {Energy
  Stable Flux Reconstruction Schemes for Advection-Diffusion Problems}},
  Computer Methods in Applied Mechanics and Engineering, 267 (2013),
  pp.~400--417, \url{https://doi.org/10.1016/j.cma.2013.08.012},
  \url{https://doi.org/10.1016/j.cma.2013.08.012}.

\bibitem{Chiew2016}
{\sc J.~J. Chiew and T.~H. Pulliam}, {\em {Stability Analysis of Dual-Time
  Stepping}}, in 46th {AIAA} Fluid Dynamics Conference, American Institute of
  Aeronautics and Astronautics, June 2016,
  \url{https://doi.org/10.2514/6.2016-3963},
  \url{https://doi.org/10.2514/6.2016-3963}.

\bibitem{Chorin1967}
{\sc A.~J. Chorin}, {\em {A Numerical Method for Solving Incompressible Viscous
  Flow Problems}}, Journal of Computational Physics, 2 (1967), pp.~12--26,
  \url{https://doi.org/10.1016/0021-9991(67)90037-x},
  \url{https://doi.org/10.1016/0021-9991(67)90037-x}.

\bibitem{Elsworth1992}
{\sc D.~Elsworth and E.~Toro}, {\em {Riemann Solvers for Solving the
  Incompressible Navier--Stokes Equations Using the Artificial Compressibility
  Method}}, Tech. Report 9208, Cranfield University, 1992.

\bibitem{Fidkowski2005}
{\sc K.~J. Fidkowski, T.~A. Oliver, J.~Lu, and D.~L. Darmofal}, {\em
  {\emph{p}-Multigrid Solution of High-Order Discontinuous Galerkin
  Discretizations of the Compressible Navier--Stokes Equations}}, Journal of
  Computational Physics, 207 (2005), pp.~92--113,
  \url{https://doi.org/10.1016/j.jcp.2005.01.005},
  \url{https://doi.org/10.1016/j.jcp.2005.01.005}.

\bibitem{Hesthaven2008}
{\sc J.~S. Hesthaven and T.~Warburton}, {\em {Nodal Discontinuous Galerkin
  Methods}}, Springer New York, 2008,
  \url{https://doi.org/10.1007/978-0-387-72067-8},
  \url{https://doi.org/10.1007/978-0-387-72067-8}.

\bibitem{Hsu2002}
{\sc J.~Hsu and A.~Jameson}, {\em {An Implicit-Explicit Hybrid Scheme for
  Calculating Complex Unsteady Flows}}, in 40th {AIAA} Aerospace Sciences
  Meeting {\&} Exhibit, American Institute of Aeronautics and Astronautics,
  Jan. 2002, \url{https://doi.org/10.2514/6.2002-714},
  \url{https://doi.org/10.2514/6.2002-714}.

\bibitem{Huynh2007}
{\sc H.~T. Huynh}, {\em {A Flux Reconstruction Approach to High-Order Schemes
  Including Discontinuous Galerkin Methods}}, in 18th {AIAA} Computational
  Fluid Dynamics Conference, American Institute of Aeronautics and
  Astronautics, June 2007, \url{https://doi.org/10.2514/6.2007-4079},
  \url{https://doi.org/10.2514/6.2007-4079}.

\bibitem{Huynh2009}
{\sc H.~T. Huynh}, {\em {A Reconstruction Approach to High-Order Schemnes
  Including Discontinuous Galerkin for Diffusion}}, in 47th {AIAA} Aerospace
  Sciences Meeting including The New Horizons Forum and Aerospace Exposition,
  American Institute of Aeronautics and Astronautics, Jan. 2009,
  \url{https://doi.org/10.2514/6.2009-403},
  \url{https://doi.org/10.2514/6.2009-403}.

\bibitem{Jameson2011}
{\sc A.~Jameson, P.~E. Vincent, and P.~Castonguay}, {\em {On the Non-linear
  Stability of Flux Reconstruction Schemes}}, Journal of Scientific Computing,
  50 (2011), pp.~434--445, \url{https://doi.org/10.1007/s10915-011-9490-6},
  \url{https://doi.org/10.1007/s10915-011-9490-6}.

\bibitem{Katzer1991}
{\sc E.~Katzer}, {\em {Multigrid Methods for Hyperbolic Equations}}, in
  Multigrid Methods {III}, Birkh\"{a}user Basel, 1991, pp.~253--263,
  \url{https://doi.org/10.1007/978-3-0348-5712-3_18},
  \url{https://doi.org/10.1007/978-3-0348-5712-3_18}.

\bibitem{Ketcheson2012}
{\sc D.~Ketcheson and A.~Ahmadia}, {\em {Optimal Stability Polynomials for
  Numerical Integration of Initial Value Problems}}, Communications in Applied
  Mathematics and Computational Science, 7 (2012), pp.~247--271,
  \url{https://doi.org/10.2140/camcos.2012.7.247},
  \url{https://doi.org/10.2140/camcos.2012.7.247}.

\bibitem{Loppi2018}
{\sc N.~Loppi, F.~Witherden, A.~Jameson, and P.~Vincent}, {\em {A High-Order
  Cross-Platform Incompressible Navier--Stokes Solver via Artificial
  Compressibility with Application to a Turbulent Jet}}, Computer Physics
  Communications, 233 (2018), pp.~193--205,
  \url{https://doi.org/10.1016/j.cpc.2018.06.016},
  \url{https://doi.org/10.1016/j.cpc.2018.06.016}.

\bibitem{Loppi2019}
{\sc N.~Loppi, F.~Witherden, A.~Jameson, and P.~Vincent}, {\em {Locally
  Adaptive Pseudo-Time Stepping for High-Order Flux Reconstruction}}, Journal
  of Computational Physics, 399 (2019), p.~108913,
  \url{https://doi.org/10.1016/j.jcp.2019.108913},
  \url{https://doi.org/10.1016/j.jcp.2019.108913}.

\bibitem{Ou2011}
{\sc K.~Ou, P.~Vincent, and A.~Jameson}, {\em {High-Order Methods for Diffusion
  Equation with Energy Stable Flux Reconstruction Scheme}}, in 49th {AIAA}
  Aerospace Sciences Meeting including the New Horizons Forum and Aerospace
  Exposition, American Institute of Aeronautics and Astronautics, Jan. 2011,
  \url{https://doi.org/10.2514/6.2011-46},
  \url{https://doi.org/10.2514/6.2011-46}.

\bibitem{Peaceman1955}
{\sc D.~W. Peaceman and J.~H.~H.~Rachford}, {\em {The Numerical Solution of
  Parabolic and Elliptic Differential Equations}}, Journal of the Society for
  Industrial and Applied Mathematics, 3 (1955), pp.~28--41,
  \url{https://doi.org/10.1137/0103003}, \url{https://doi.org/10.1137/0103003}.

\bibitem{Peyret1976}
{\sc R.~Peyret}, {\em {Unsteady Evolution of a Horizontal Jet in a Stratified
  Fluid}}, Journal of Fluid Mechanics, 78 (1976), pp.~49--63,
  \url{https://doi.org/10.1017/s0022112076002322},
  \url{https://doi.org/10.1017/s0022112076002322}.

\bibitem{Rogers1995}
{\sc S.~Rogers}, {\em {A Comparison of Implicit Schemes for the Incompressible
  Navier-Stokes Equations with Artificial Compressibility}}, in 33rd Aerospace
  Sciences Meeting and Exhibit, American Institute of Aeronautics and
  Astronautics, Jan. 1995, \url{https://doi.org/10.2514/6.1995-567},
  \url{https://doi.org/10.2514/6.1995-567}.

\bibitem{Toro2009}
{\sc E.~F. Toro}, {\em {Riemann Solvers and Numerical Methods for Fluid
  Dynamics}}, Springer Berlin Heidelberg, 2009,
  \url{https://doi.org/10.1007/b79761}, \url{https://doi.org/10.1007/b79761}.

\bibitem{Trojak2018}
{\sc W.~Trojak, R.~Watson, and P.~G. Tucker}, {\em {Temporal Stabilisation of
  Flux Reconstruction on Linear Problems}}, in 2018 Fluid Dynamics Conference,
  American Institute of Aeronautics and Astronautics, June 2018,
  \url{https://doi.org/10.2514/6.2018-4263},
  \url{https://doi.org/10.2514/6.2018-4263}.

\bibitem{Vermeire2019}
{\sc B.~Vermeire, N.~Loppi, and P.~Vincent}, {\em {Optimal Runge--Kutta Schemes
  for Pseudo Time-Stepping with High-Order Unstructured Methods}}, Journal of
  Computational Physics, 383 (2019), pp.~55--71,
  \url{https://doi.org/10.1016/j.jcp.2019.01.003},
  \url{https://doi.org/10.1016/j.jcp.2019.01.003}.

\bibitem{Vincent2010}
{\sc P.~E. Vincent, P.~Castonguay, and A.~Jameson}, {\em {A New Class of
  High-Order Energy Stable Flux Reconstruction Schemes}}, Journal of Scientific
  Computing, 47 (2010), pp.~50--72,
  \url{https://doi.org/10.1007/s10915-010-9420-z},
  \url{https://doi.org/10.1007/s10915-010-9420-z}.

\bibitem{Wienands2001}
{\sc R.~Wienands and C.~W. Oosterlee}, {\em {On Three-Grid Fourier Analysis for
  Multigrid}}, {SIAM} Journal on Scientific Computing, 23 (2001), pp.~651--671,
  \url{https://doi.org/10.1137/s106482750037367x},
  \url{https://doi.org/10.1137/s106482750037367x}.

\bibitem{Williams2014}
{\sc D.~Williams, L.~Shunn, and A.~Jameson}, {\em {Symmetric Quadrature Rules
  for Simplexes Based on Sphere Close Packed Lattice Arrangements}}, Journal of
  Computational and Applied Mathematics, 266 (2014), pp.~18--38,
  \url{https://doi.org/10.1016/j.cam.2014.01.007},
  \url{https://doi.org/10.1016/j.cam.2014.01.007}.

\bibitem{Witherden2014}
{\sc F.~Witherden, A.~Farrington, and P.~Vincent}, {\em {{Pyfr}: An Open Source
  Framework for Solving Advection-Diffusion Type Problems on Streaming
  Architectures Using the Flux Reconstruction Approach}}, Computer Physics
  Communications, 185 (2014), pp.~3028--3040,
  \url{https://doi.org/10.1016/j.cpc.2014.07.011},
  \url{https://doi.org/10.1016/j.cpc.2014.07.011}.

\bibitem{Yoon1987}
{\sc S.~Yoon and A.~Jameson}, {\em {An {LU}-{SSOR} Scheme for the Euler and
  Navier--Stokes Equations}}, in 25th {AIAA} Aerospace Sciences Meeting,
  American Institute of Aeronautics and Astronautics, Mar. 1987,
  \url{https://doi.org/10.2514/6.1987-600},
  \url{https://doi.org/10.2514/6.1987-600}.

\end{thebibliography}


\begin{appendices}
\section{FR Operator Definition}
\label{app:fr_op}
    The FR operators of first-order derivatives are defined as 
    \begin{subequations}
        \begin{align}
            \px{\mathbf{u}_i}{x} &= \frac{2}{h}\left(\mathbf{C}_-\mathbf{u}_{i-1} + \mathbf{C}_0\mathbf{u}_i + \mathbf{C}_+\mathbf{u}_{i+1}\right),\\
            \mathbf{C}_- &= \alpha\mathbf{g}_L\mathbf{l}_R^T,\\
            \mathbf{C}_0 &= \mathbf{D} - \alpha\mathbf{g}_L\mathbf{l}_L^T - (1-\alpha)\mathbf{g}_R\mathbf{l}_R^T,\\
            \mathbf{C}_+ &= \alpha\mathbf{g}_L\mathbf{l}_R^T.
        \end{align}
    \end{subequations}
    The matrix $\mathbf{D}$ is the nodal differentiation matrix ($\mathbf{D}_{ij}=\pxvar{l_i(x_j)}{x}$), $\mathbf{g}_L$ is the \emph{gradient} of the left correction function at the solution points and $\mathbf{l}_l$ is the interpolation of the solution points to the left faces.
    
    The FR methodology for second-order derivatives is to nest the derivatives, treating each first order derivative in the standard manner~\cite{Castonguay2013,Huynh2009,Ou2011}. In particular, the diffusion equation takes the form
    \begin{equation}
        \px{u}{t} = \mu\px{q}{x}, \quad \mathrm{where} \quad  q = \px{u}{x}.
    \end{equation}
    Each stage is then solved with the FR methodology, which in the vector form is
    \begin{subequations}
        \begin{align}
            \mathbf{q}_i &= \frac{2}{h}\Big(\mathbf{C}_{-}\mathbf{u}_{i-1} + \mathbf{C}_0\mathbf{u}_i + \mathbf{C}_{+}\mathbf{u}_{i+1}\Big),\\
            \px{\mathbf{q}_i}{x} &= \frac{2}{h}\Big(\mathbf{C}_{-}\mathbf{q}_{i-1} + \mathbf{C}_0\mathbf{q}_i + \mathbf{C}_{+}\mathbf{q}_{i+1}\Big).
        \end{align}
    \end{subequations}
    These may be combined to achieve
    \begin{multline}
        \pxi{2}{\mathbf{u}_i}{x} = \mathbf{Q}_d\mathbf{u}_i = \frac{4}{h^2}\Big(\mathbf{C}_{-}^2\mathbf{u}_{i-2} + 
		    (\mathbf{C}_{-}\mathbf{C}_0 + \mathbf{C}_0\mathbf{C}_{-})\mathbf{u}_{i-1} + \\
		    (\mathbf{C}_{-}\mathbf{C}_{+} + \mathbf{C}_0^2 + \mathbf{C}_{+}\mathbf{C}_{-})\mathbf{u}_{i} + \\
		    (\mathbf{C}_0\mathbf{C}_{+} + \mathbf{C}_{+}\mathbf{C}_0)\mathbf{u}_{i+1} + 
		    \mathbf{C}_{+}^2\mathbf{u}_{i+2}\Big).
    \end{multline}
    In the analysis performed in the main body of this work the following assignments are used for brevity.
    \begin{subequations}
        \begin{align}
            \mathbf{B}_{-2} &= \mathbf{C}_{-}^2\\
            \mathbf{B}_{-} &= \mathbf{C}_{-}\mathbf{C}_0 + \mathbf{C}_0\mathbf{C}_{-}\\
            \mathbf{B}_{0} &= \mathbf{C}_{-}\mathbf{C}_{+} + \mathbf{C}_0^2 + \mathbf{C}_{+}\mathbf{C}_{-}\\
            \mathbf{B}_{+} &= \mathbf{C}_0\mathbf{C}_{+} + \mathbf{C}_{+}\mathbf{C}_0\\
            \mathbf{B}_{+2} &= \mathbf{C}_{+}^2
        \end{align}
    \end{subequations}
\section{Element Jacobi Smoothing}
\label{app:ej}
    Fidkowski~\etal~\cite{Fidkowski2005} investigated the use of \emph{p}-multigrid on the convergence of implicit DG with Element-Jacobi (EJ) smoothing. As a canonical approach for solving implicit systems of equations we have included this method to provide a benchmark for the dual time approach. The equivalent of the pseudo time update for EJ takes the form
    \begin{subequations}
        \begin{align}
            \mathbf{u}_{n+1,m+1} &= \mathbf{u}_{n+1,m} - \kappa\mathbf{J}^{-1}(\mathbf{T}\mathbf{u}_{n+1,m} - C_B\mathbf{u}_{n+1,0}),\\
            \mathbf{J} &= \px{}{\mathbf{u}_{n+1,m}}(\mathbf{T}\mathbf{u}_{n+1,m} - C_B\mathbf{u}_{n+1,0}),
        \end{align}
    \end{subequations}
    where $\kappa$ is the relaxation factor. From \cref{eq:Q_def,eq:pseudo_res} the Jacobian matrix inverse may then be defined as
    \begin{equation}
        \mathbf{J}^{-1} = B_0\Delta t\left[\mathbf{I} - \frac{2B_0\Delta t}{h}\left(\mathbf{C}_0 -\frac{2\mu}{h}\mathbf{B}_0\right)\right]^{-1},
    \end{equation}
    and this may then be inserted in the previously defined \emph{p}-multigrid algorithms in place of the RK pseudo-time integration.
\end{appendices}


\end{document}